\documentclass{article}
\usepackage{a4wide}
\usepackage{graphicx}
\usepackage{amssymb}
\usepackage{amsmath}

\providecommand{\tabularnewline}{\\}

 \usepackage{amsthm}
 \theoremstyle{plain}
 \newtheorem{thm}{Theorem}[section]
 \numberwithin{equation}{section}
 \numberwithin{figure}{section}
 \theoremstyle{plain}
 \newtheorem{lem}[thm]{Lemma} 
 \theoremstyle{definition}
 \newtheorem*{terminology}{Terminology}
 \theoremstyle{plain}
 \newtheorem{cor}[thm]{Corollary}
 \theoremstyle{plain}

\newcommand{\OEDc}{\mathcal{E}_\mathrm{o}^\mathrm{o}}
\newcommand{\REDc}{\mathcal{E}_\mathrm{o}^\mathrm{r}}
\newcommand{\IDc}{\mathcal{E}_\mathrm{i}}
\newcommand{\SDc}{\mathcal{E}_\mathrm{s}}

\newcommand{\VDc}{\mathcal{V}}
\newcommand{\OFDc}{\mathcal{F}^\mathrm{o}}
\newcommand{\FDc}{\mathcal{F}}
\newcommand{\SFDc}{\mathcal{F}_\mathrm{s}}
\newcommand{\Dc}{\mathcal{D}}
\newcommand{\VCc}{\hat{\mathcal{C}}}
\newcommand{\Cc}{\mathcal{C}}
\newcommand{\Outc}{\mathcal{G}}

\newcommand{\OED}{E_\mathrm{o}^\mathrm{o}}
\newcommand{\RED}{E_\mathrm{o}^\mathrm{r}}
\newcommand{\ID}{E_\mathrm{i}}
\newcommand{\SD}{E_\mathrm{s}}

\newcommand{\VD}{V}
\newcommand{\OFD}{F^\mathrm{o}}
\newcommand{\FD}{F}
\newcommand{\SFD}{F_\mathrm{s}}
\newcommand{\D}{D}
\newcommand{\VC}{\hat{C}}
\newcommand{\C}{C}
\newcommand{\Out}{G}

\newcommand{\vd}{v}

\newcommand{\dis}{d}
\newcommand{\vc}{\hat{c}}
\newcommand{\con}{c}
\newcommand{\out}{g}

\newcommand{\bCc}{\mathcal{C}_\mathrm{\it{b}}}

\newcommand{\bC}{C_\mathrm{\it{b}}}

\newcommand{\bout}{g_\mathrm{\it{b}}}

\begin{document}

\title{Enumeration of Unlabeled Outerplanar Graphs}

\author{\vspace*{2em}Manuel Bodirsky$^{1}$, \'Eric Fusy$^{2}$, Mihyun Kang$^{1,3}$, and Stefan Vigerske$^{1}$
 \\
 $^{1}$Humboldt-Universit\"at zu Berlin, Institut f\"ur Informatik\\
 Unter den Linden 6, 10099 Berlin, Germany\\
 \{bodirsky, kang\}@informatik.hu-berlin.de\\
 stefan@mathematik.hu-berlin.de\vspace*{1em}\\
 $^{2}$Projet Algo, INRIA Rocquencourt\\
 B. P. 105, 78153 Le Chesnay Cedex, France\\
 eric.fusy@inria.fr}

\maketitle

\footnotetext[3]{Research supported by the Deutsche Forschungsgemeinschaft (DFG Pr 296/7-3)}

\begin{abstract}
We determine the exact and asymptotic
number of \emph{unlabeled outerplanar graphs}.
The exact number $\out_{n}$ of unlabeled outerplanar graphs 
on $n$ vertices can be computed in 
polynomial time, and $\out_{n}$ is asymptotically $\out\,
n^{-5/2}\rho^{-n}$, where $\out\approx0.00909941$ and
$\rho^{-1}\approx7.50360$ can be approximated.
Using our enumerative results we investigate several statistical properties of
 random unlabeled outerplanar graphs on $n$ vertices, for instance 
concerning connectedness, chromatic number, and the number of edges.
To obtain the results we combine classical cycle index enumeration with
recent results from analytic combinatorics.\\
\begin{small}\textit{Keywords:} 
unlabeled outerplanar graphs, dissections, combinatorial enumeration, cycle index, asymptotic estimates, singularity analysis
\end{small}
\end{abstract}

\section{Introduction and results}
Singularity analysis is a highly successful tool for asymptotic
enumeration of combinatorial structures~\cite{FlajoletSedgewick}, 
once we have a sufficiently good description of the corresponding 
exponential or ordinary generating functions.
If we want to count \emph{unlabeled} structures, i.e., if we count the 
structures up to isomorphism, the potential symmetries of the structures
often require a more powerful tool than generating functions, 
namely \emph{cycle indices}, introduced by P\'olya~\cite{Polya}.
From the cycle index sums for a class of combinatorial structures
we can obtain its generating function, and apply singularity analysis.
However, when the cycle index sums are given only implicitly, it might be a challenging task to apply this technique.
This is illustrated by the situation for planar graphs: 
the asymptotic number of labeled planar graphs
was recently determined by Gim\'enez and Noy~\cite{GimenezNoy}, 
based on singularity analysis, whereas the
enumeration of unlabeled planar graphs is a research problem left
open for several decades~\cite{tensteps}. 

In this paper we
determine the exact and asymptotic
number of \emph{unlabeled outerplanar graphs}, 
an important subclass of the class
of all planar graphs. 
We provide a polynomial-time algorithm to compute the exact number
$\out_{n}$ of unlabeled outerplanar graphs on $n$ vertices, and 
prove that $\out_{n}$ is asymptotically $\out\,
n^{-5/2}\rho^{-n}$, where $\out\approx0.00909941$ and
$\rho^{-1}\approx7.50360$ can be approximated.
Building on our enumerative results we derive typical properties of 
 random unlabeled outerplanar graphs on $n$ vertices (i.e. taken uniformly at random among all unlabeled outerplanar graphs on $n$ vertices), for example connectedness, chromatic number, the number of components, and the number of edges.

Before we provide a more detailed exposition of the results obtained in 
this paper, we would like to give a brief survey on the vast literature
on enumerative results for planar structures. 
The exact and asymptotic number of \emph{embedded} planar graphs (i.e., planar maps) 
has been studied intensively,
starting with Tutte's seminal work on the number of rooted oriented planar maps~\cite{Tutte}. The number of three-connected planar maps is related to the number of three-connected planar graphs~\cite{MullinSchellenberg,Tutte}, since a three-connected planar graph has a unique embedding on the sphere~\cite{Whitney}. 
Bender, Gao, and Wormald used this property 
to count labeled two-connected planar graphs \cite{BenderGaoWormald}, 
and Gim\'enez and Noy recently extended this work to the enumeration 
of labeled planar graphs~\cite{GimenezNoy}. 
The growth constant for labeled planar graphs can be computed 
with arbitrary precision and its first digits are $27.2269$. 
Many interesting properties of a random labeled planar graph were studied by McDiarmid, Steger, and Welsh~\cite{McDiarmid}. 
It is also known how to generate labeled three-connected planar graphs,  labeled planar maps, and labeled planar graphs uniformly at random \cite{BodirskyGroepletal,BodirskyGroepletal4,FPS,Schaeffer}.

The asymptotic number of general unlabeled planar graphs 
has not yet been determined, but has been studied 
for quite some time~\cite{tensteps}.
Moreover, no polynomial time algorithm for the computation of the exact 
number of unlabeled planar graphs on $n$ vertices 
is known. Such an algorithm is only
known for unlabeled rooted two-connected planar graphs~\cite{BodirskyGroepletal2}, and for 
unlabeled rooted 
cubic planar graphs~\cite{BodirskyGroepletal3}. 

An \emph{outerplanar graph} is a graph that can be embedded in the plane such
that every vertex lies on the outer face. Such graphs can also be
characterized in terms of forbidden minors~\cite{ChartrandHarary}, namely $K_{2,3}$ and $K_{4}$. The class of outerplanar graphs is often used as a first
non-trivial test-case for results about the class of all planar graphs; 
apart from that, 
this class appears frequently in various applications of graph theory.
The asymptotic number of \emph{labeled} outerplanar graphs 
was recently determined in~\cite{BodirskyGimenezKangNoy}. 
In this paper, we determine the number of \emph{unlabeled} outerplanar graphs,
i.e., we enumerate outerplanar graphs up to isomorphism. 

Two-connected outerplanar graphs can
be identified with dissections of a convex
polygon~\cite{Whitney}. General outerplanar graphs can be decomposed 
according to
their degree of connectivity: an outerplanar
graph is a set of connected outerplanar graphs, and a connected
outerplanar graph can be decomposed into two-connected blocks. In the
labeled case this
decomposition yields equations that link the exponential generating 
functions of
two-connected, connected, and general outerplanar 
graphs~\cite{BodirskyGimenezKangNoy}. Once labeled
dissections are enumerated, these equations yield formulas
for counting outerplanar graphs.  
In the unlabeled case we use the same decomposition. 
However, the potential symmetries make
it more difficult to obtain exact and asymptotic results.
We have to use cycle index sums, which were introduced by P\'olya for 
unlabeled enumeration~\cite{Polya}, to obtain implicit
information about the ordinary generating functions of unlabeled
outerplanar graphs. We then apply
singularity analysis, a very powerful tool that is thoroughly developed in the forthcoming book of Flajolet and Sedgewick~\cite{FlajoletSedgewick}.
A similar strategy was applied by Labelle, Lamathe, and Leroux for the
enumeration of unlabeled $k$-gonal 2-trees~\cite{LLL}.
However, in the singularity analysis for outerplanar graphs we face 
new difficulties, which did not appear in other literature, as far as we
know. The generating function for connected outerplanar graphs is
defined implicitly by a multiset of connected outerplanar graphs that
are rooted at a two-connected component, and moreover the number of
two-connected graphs has exponential growth.  By applying the singular
implicit function theorem we overcome these difficulties (see Section~\ref{sec:estimates} for the details) and estimate the asymptotic number
of outerplanar graphs.

\paragraph*{Contributions.}
Our first result is the extension of Read's counting formulas~\cite{Read} 
for the number of unlabeled two-connected outerplanar graphs 
to counting formulas for the number of unlabeled outerplanar graphs.

\begin{thm}
\label{thm:polyalg}
The exact numbers of unlabeled two-connected outerplanar graphs $\dis_{n}$, unlabeled connected outerplanar graphs $\con_{n}$, and unlabeled outerplanar graphs $\out_{n}$ with $n$ vertices can be computed in polynomial time.
\end{thm}
See the sequences A001004, A111563, and A111564 from~\cite{Sloane} for initial values. 

\begin{thm}\label{thm:asymptoticouter}
The numbers $\dis_{n}$, $\con_{n}$, and $\out_{n}$ of two-connected, connected, and general outerplanar graphs with $n$ vertices have the asymptotic estimates
\begin{eqnarray*}
\dis_{n} & \sim & \dis\, n^{-\frac{5}{2}}\delta^{-n},\\
\con_{n} & \sim & \con\, n^{-\frac{5}{2}}\rho^{-n},\\
\out_{n} & \sim & \out\, n^{-\frac{5}{2}}\rho^{-n},
\end{eqnarray*}
with growth rates $\delta^{-1}=3+2\sqrt{2}\approx5.82843$ and $\rho^{-1}\approx7.50360$,
and constants $\dis \approx0.00596026$, $\con\approx0.00760471$,
and $\out\approx0.00909941$. (See Theorems \ref{thm:twoconn}, \ref{thm:Casympt},
and \ref{thm:Out_asympt}.)
\end{thm}

Having the asymptotic estimates of unlabeled connected outerplanar
graphs and unlabeled outerplanar graphs, we investigate  asymptotic
distributions of parameters such as the number of components and the
number of isolated vertices of a random outerplanar graph on $n$ vertices.
\begin{thm}
\label{thm:typprop}
\begin{enumerate}
\item[(1)] The probability that a random outerplanar graph is connected is
asymptotically $\con/\out\approx0.845721$.
\item[(2)] The expected number of components in a random unlabeled outerplanar graph is asymptotically equal to a constant $\approx1.17847$.
\item[(3)] The asymptotic distribution of the number of isolated
vertices in a random outerplanar graph is a geometric law with parameter $\rho$. In particular, the expected number of isolated vertices in a random outerplanar graph is asymptotically $\rho/\left(1-\rho\right)\approx0.153761$.
\end{enumerate}
\end{thm}

To investigate the chromatic number of a random outerplanar graph we also study  the asymptotic number of unlabeled bipartite outerplanar graphs.

\begin{thm}
\label{thm:bdisgrowthrate}
The number of bipartite outerplanar graphs $\left(\bout\right)_n$ on $n$ vertices has the asymptotic estimate
$\left(\bout\right)_{n} \sim b n^{-5/2}\rho_{\mathrm{\it{b}}}^{-n}$,
with $\rho_{\mathrm{\it{b}}}^{-1}\approx 4.57717$.
\end{thm}

An outerplanar graph is easily shown to have a 3-colouring. The fact that the growth constant of bipartite outerplanar graphs is smaller than the growth constant of outerplanar graphs yields the following result:
\begin{thm}
%[chromatic number]
\label{thm:chromaticnumber}
The probability that the chromatic number of a random
unlabeled outerplanar graph is different from three decays asymptotically exponentially to zero.
\end{thm}

If we count graphs with respect to the additional parameter that specifies the number of edges, we can study the distribution of the number of edges in a random outerplanar graph.
\begin{thm}
\label{thm:edges}
The distribution
of the number of edges in a random outerplanar graph on $n$ vertices
is asymptotically Gaussian with mean $\mu n$ and variance $\sigma^2n$, where $\mu\approx1.54894$ and $\sigma^2\approx0.227504$. The same holds for random connected outerplanar
graphs with the same mean and variance and for
random two-connected outerplanar graphs with asymptotic mean $\left(1+\sqrt{2}/2\right)n\approx1.70711n$
and asymptotic variance $\sqrt{2}/8\, n\approx0.176777n$.
\end{thm}

\paragraph*{Outline.}
The paper is organized as follows. Section \ref{sec:preliminaries} introduces well-known techniques for
the enumeration of rooted and unrooted unlabeled structures and shows how to obtain asymptotic estimates.
Section \ref{sec:cycleindexsums} provides exact enumeration of two-connected, connected, and general outerplanar graphs, and also of bipartite outerplanar graphs.
Section~\ref{sec:asymptotics}
provides asymptotic estimates for the number of two-connected, connected, general,
and bipartite outerplanar graphs. Section~\ref{sec:approxgrowth} shows how to approximate the growth
constant for outerplanar graphs. Finally, Section~\ref{sec:randomgraphs} investigates typical properties
of a random outerplanar graph on $n$ vertices, such as the probability of connectedness, the expected
number of components, the expected number of isolated vertices, the chromatic number, and the distribution of the number of edges.

\section{\label{sec:preliminaries}Preliminaries}
We recall some concepts and techniques that we need for the enumeration of unlabeled graphs, and some facts from singularity analysis to obtain asymptotic estimates.

\subsection{Cycle index sums}
To enumerate unlabeled graphs, we use \emph{cycle indices} as introduced by P\'olya \cite{HararyPalmer,Polya}. For a group of permutations $A$ on an object set $X=\left\{ 1,\ldots,n\right\} $ (for example, the vertex set of a graph), the cycle index $Z\left(A\right)$ of $A$ with respect to the formal variables $s_{1},\ldots,s_{n}$ is defined by
\[
Z\left(A\right):=Z\left(A;s_{1},s_{2},\ldots\right):=\frac{1}{\left|A\right|}\sum_{\alpha\in A}\prod_{k=1}^{n}s_{k}^{j_{k}\left(\alpha\right)},\]
where $j_{k}\left(\alpha\right)$ denotes the number of cycles of length $k$ in the disjoint cycle decomposition of $\alpha\in A$. For a graph $G$ on $n$ vertices with automorphism group $\Gamma\left(G\right)$, we write $Z\left(G\right):=Z\left(\Gamma\left(G\right)\right)$, and
for a set of graphs $\mathcal{K}$, we write $Z\left(\mathcal{K}\right)$ for the \emph{cycle index sum} for $\mathcal{K}$ defined by
\[Z\left(\mathcal{K}\right):=Z\left(\mathcal{K};s_{1},s_{2},\ldots\right):=\sum_{K\in\mathcal{K}}Z\left(K;s_{1},s_{2},\ldots\right).\]
It can be shown~\cite{BergeronLabelleLeroux} that, if $\bar{\mathcal{K}}$ is the set of graphs of $\mathcal{K}$ equipped with distinct labels, then
\[Z(\mathcal{K})=\sum_{n\geq 0}\frac{1}{n!}\sum_{K\in \bar{\mathcal{K}}_n}\sum_{\alpha\in\Gamma(K)}\prod_{k=1}^n s_k^{j_k(\alpha)},
\]
which coincides with the classical definition of a \emph{cycle index series} and shows the close relationship
of cycle index sums to exponential generating functions in labeled counting.

Indeed, cycle index sums can be used for the enumeration
of unlabeled structures in a similar way as generating functions for labeled enumeration.
First of all, the composition of graphs corresponds to the composition of the associated cycle indices. 
Consider an object set $X=\left\{ 1,\ldots,n\right\} $ and a permutation group $A$ on $X$. A composition of $n$ graphs from $\mathcal{K}$ is a function $f:X\rightarrow\mathcal{K}$. Two compositions $f$ and $g$ are similar, $f\sim g$, if there exists a permutation $\alpha\in A$ with $f\circ\alpha=g$. We write $\mathcal{G}$ for the set of equivalence classes of compositions of $n$ graphs from $\mathcal{K}$ (with respect to the equivalence relation $\sim$). Then
\begin{equation}
Z\left(\mathcal{G}\right)=Z\left(A\right)\left[Z\left(\mathcal{K}\right)\right]:=Z\left(A;Z\left(\mathcal{K};s_{1},s_{2},\ldots\right),Z\left(\mathcal{K};s_{2},s_{4},\ldots\right),\ldots\right),\label{eq:composition}
\end{equation}
i.e., $Z(\mathcal{G})$ is obtained from $Z(A)$ by replacing each $s_{i}$ by $Z\left(\mathcal{K};s_{i},s_{2i},\ldots\right)=\sum_{K\in\mathcal{K}}Z\left(K;s_{i},s_{2i},\ldots\right)$ \cite{HararyPalmer}. Hence, Formula (\ref{eq:composition}) makes it possible to derive the cycle index sum for a class of graphs by decomposing the graphs into simpler structures with known cycle index sum.

In many cases, such a decomposition is only possible when, for example, one vertex is distinguished from the others in the graphs, so that there is a unique point where the decomposition is applied. Graphs with a distinguished vertex are called vertex \emph{rooted} graphs. The automorphism group of a vertex rooted graph consists of all permutations of the group of the unrooted graph that fix the root vertex. Hence, one can expect a close relation between the cycle index of unrooted graphs and the cycle indices of their rooted counterparts. As shown in \cite{HararyPalmer}, if $\mathcal{G}$ is an unlabeled set of graphs and $\hat{\mathcal{G}}$ is the set of graphs of $\mathcal{G}$ rooted at a vertex, then
\begin{equation}
Z(\hat{\mathcal{G}})=s_{1}\frac{\partial}{\partial s_{1}}Z\left(\mathcal{G}\right).\label{eq:rooting}
\end{equation}
This relationship can be inverted to express the cycle index sum for
the unrooted graphs in terms of the cycle index sum for the rooted graphs,
\begin{equation}
Z\left(\mathcal{G}\right)=\int_{0}^{s_{1}}\frac{1}{s_{1}}Z(\hat{\mathcal{G}})ds_{1}+Z\left(\mathcal{G}\right)\left|_{s_{1}=0}\right..\label{eq:unrooting}\end{equation}
Observe that permutations without fixed points are not counted by the cycle indices of the rooted graphs, so that their cycle indices are added as a boundary term to $Z\left(\mathcal{G}\right)$.

\subsection{\label{sub:ogf}Ordinary generating functions}
Once the cycle index sum for a class of graphs of interest is known, the corresponding ordinary generating function can be derived by replacing the formal variables $s_{i}$ in the cycle index sums by $x^{i}$ (note that $Z\left(G;x,x^{2},\ldots\right)=x^{\left|G\right|}$ for a graph $G$). More generally, for a group $A$ and an ordinary generating function $K(x)$ we define
\[
Z\left(A;K(x)\right):=Z\left(A;K(x),K(x^{2}),K(x^{3}),\ldots\right)\]
as the ordinary generating function obtained by substituting each $s_{i}$ in $Z\left(A\right)$ by $K(x^{i})$, $i\geq1$.

\subsection{The dissimilarity characteristic theorem}
The dissimilarity characteristic theorem expresses the number of dissimilar vertices of a graph in terms of the numbers of dissimilar blocks and the number of dissimilar vertices of each block in the graph \cite{HararyPalmer}. In the case of trees, the blocks of the graph are the edges. An edge whose vertices are interchanged by an automorphism of the graph is called a \textit{symmetry-edge}, and has to be treated separately. The resulting equation, which relates the number of dissimilar vertices to the number of dissimilar edges and symmetry-edges in a tree, can be used to derive the generating function for unrooted trees from the generating functions for trees rooted at a vertex, at an edge, or at a symmetry-edge, respectively \cite{HararyPalmer}. To also obtain the cycle index sum for unrooted trees, the dissimilarity characteristic theorem can be extended by considering vertex rooted trees whose root vertex is incident
to a symmetry-edge \cite{Vigerske}.

\begin{lem}
\label{lem:dissimilar_cycleindex}
Let $G$ be an unlabeled tree. Let $\mathcal{V}$ be the set of vertex rooted trees
that have $G$ as underlying unrooted tree. Partition the set $\mathcal{V}$ 
into the set $\mathcal{V}_{1}$ of rooted trees  where the root vertex is not incident to a symmetry-edge and the set $\mathcal{V}_{2}$ where the root vertex is incident to a symmetry-edge. Furthermore, let $\mathcal{E}$ be the set of trees obtained from $G$ by rooting at an edge, and let $\mathcal{S}$ be the set of trees obtained from $G$ by rooting at a symmetry-edge. Then
\begin{equation}
Z\left(G\right)=Z\left(\mathcal{V}_{1}\right)-Z\left(\mathcal{E}\right)+2Z\left(\mathcal{S}\right).\label{eq:cycleindex_disstheorem}
\end{equation}
\end{lem}
The lemma can be proven by induction on the number of dissimilar edges, where the initial case of the induction determines the crucial distinction between the two kinds of vertex rooted trees, $\mathcal{V}_{1}$ and $\mathcal{V}_{2}$ (see \cite{Vigerske} for a full proof).

\subsection{\label{sub:singularityanalysis}Singularity analysis}
To determine asymptotic estimates of the coefficients of a generating function we use singularity analysis~\cite{FlajoletSedgewick}. The fundamental observation is that the exponential growth of the coefficients of a function that is analytic at the origin is determined by the dominant singularities of the function, i.e., singularities at the boundary of the disc of convergence. By Pringsheim's theorem \cite[Thm. IV.6]{FlajoletSedgewick},
a generating function $F\left(z\right)$ with non-negative coefficients and finite radius of convergence $R$ has a singularity at the point $z=R$. If $z=R$ is the unique singularity on the disk $|z|=R$, it follows from the exponential growth formula \cite[Thm. IV.7]{FlajoletSedgewick} that the coefficients $f_{n}=\left[z^{n}\right]F\left(z\right)$ satisfy $f_{n}=\theta\left(n\right)R^{-n}$ with $\limsup_{n\rightarrow\infty}\left|\theta\left(n\right)\right|^{1/n}=1$. A closer look at the type of the dominant singularity, for example, the order of the pole, enables the computation of subexponential factors as well. The following lemma describes the singular expansion for a common case \cite[Thm. VI.1]{FlajoletSedgewick}.

\begin{lem}[standard function scale]\label{lem:basicscale}
Let $F\left(z\right)=\left(1-z\right)^{-\alpha}$ with $\alpha\not\in\left\{ 0,-1,-2,\ldots\right\} $. Then the coefficients $f_{n}$ of $F(z)$ have a full asymptotic development in descending powers of $n$,
\begin{equation}
f_{n}=\binom{n+\alpha-1}{n}
\sim\frac{n^{\alpha-1}}{\Gamma\left(\alpha\right)}\left(1+\sum_{k=1}^{\infty}\frac{e_{k}\left(\alpha\right)}{n^{k}}\right)\label{eq:basisscale}
\end{equation}
where $\Gamma\left(\alpha\right)$ is the Gamma-Function, $\Gamma\left(\alpha\right):=\int_{0}^{\infty}e^{-t}t^{\alpha-1}dt$ for $\alpha\not\in\left\{ 0,-1,-2,\ldots\right\} $, and $e_{k}\left(\alpha\right)$ is a polynomial in $\alpha$ of degree $2k$.
\end{lem}
In our calculations, it will appear that a generating function $f(x)$ is given only implicitly by an equation $H(x,f(x))=0$.  Theorem VII.3 in \cite{FlajoletSedgewick} describes how to derive a full singular expansion of $f(x)$ in this case. We state it here in a slightly modified version. A generating function is called \emph{aperiodic}, if it can not be written in the form $Y\left(x\right)=x^{a}\tilde{Y}\left(x^{d}\right)$ with $d\geq2$ and $\tilde{Y}$ analytic at $0$.

\begin{thm}
[singular implicit functions]
\label{thm:positiveimplicitfunctions}Let
$H\left(x,y\right)$ be a bivariate function that is analytic in a complex domain $|x|<R$, $|y|<S$ and verifies $H(0,0)=0$, $\frac{\partial}{\partial y}H\left(0,0\right)=-1$, and whose Taylor coefficients $h_{m,n}$ satisfy the following positivity conditions: they are nonnegative except for $h_{0,1}=-1$ (because $\frac{\partial}{\partial y}H\left(0,0\right)=-1$) and $h_{m,n}>0$ for at least one pair $(m,n)$ with $n \geq 2$.
 Assume that there are two numbers
$r\in\left(0,R\right)$ and $s\in\left(0,S\right)$ such that\begin{equation}
H\left(r,s\right)=0,\qquad\frac{\partial}{\partial y}H\left(r,s\right)=0,\label{eq:posimplicitfunc}\end{equation}
$\frac{\partial^{2}}{\partial y^{2}}H\left(r,s\right)\neq0$ and $\frac{\partial}{\partial x}H\left(r,s\right)\neq0$.
Assume further that the equation $H\left(x,Y\left(x\right)\right)=0$
admits a solution $Y\left(x\right)$ that is analytic at $0$, has
non-negative coefficients, and is aperiodic.
Then $r$ is the unique dominant singularity of $Y(x)$ and $Y\left(x\right)$ converges at $x=r$, where
it has the singular expansion\[
Y\left(x\right)=s+\sum_{i\geq1}Y_{i}\left(\sqrt{1-\frac{x}{r}}\right)^{i},\qquad\textrm{with }Y_{1}=-\sqrt{\frac{2r\frac{\partial}{\partial x}H\left(r,s\right)}{\frac{\partial^{2}}{\partial y^{2}}H\left(r,s\right)}}\neq0,\]
and computable constants $Y_{2},Y_{3},\cdots$.
Hence,\[\left[x^{n}\right]Y\left(x\right)=-\frac{Y_{1}}{2\sqrt{\pi n^{3}}}r^{-n}\left(1+O\left(\frac{1}{n}\right)\right).\]

\end{thm}
The formulas that express the coefficients $Y_{i}$ in terms of partial
derivatives of $H\left(x,y\right)$ at $\left(r,s\right)$ can be
found in \cite{Finch2,PlotkinRosenthal}.

When a parameter $\xi$ of a combinatorial structure is studied, the generating function $F\left(x\right)$ has to be extended to a bivariate
generating function $F\left(x,y\right)=\sum_{n,m}f_{n,m}x^{n}y^{m}$
where the second variable $y$ marks $\xi$.
We can determine the asymptotic distribution of $\xi$ from $F\left(x,y\right)$ by varying $y$
in some neighbourhood of $1$. The following theorem follows
from the so-called quasi-powers theorem~\cite[Thm. IX.7]{FlajoletSedgewick}.

\begin{thm}
\label{thm:centrallaw}
Let $F\left(x,y\right)$ be a bivariate generating
function of a family of objects $\mathcal{F}$, where the power in $y$
corresponds to a parameter $\xi$
on $\mathcal{F}$, i.e., $\left[x^{n}y^{m}\right]F\left(x,y\right)=\left|\left\{ F\in\mathcal{F}\left|\left|F\right|=n,\xi\left(F\right)=m\right.\right\} \right|$.
Assume that, in a fixed complex neighbourhood of $y=1$,  $F\left(x,y\right)$
has a singular expansion of the form\begin{equation}
F\left(x,y\right)=\sum_{k\geq0}F_{k}\left(y\right)\left(\sqrt{1-\frac{x}{x_{0}\left(y\right)}}\right)^{k}\label{eq:Fsingexpand}\end{equation}
 where $x_{0}\left(y\right)$ 
is the dominant singularity of $x\mapsto F\left(x,y\right)$.
Furthermore, assume that there is an odd $k_{0}\in\mathbb{N}$ 
such that for all $y$ in the neighbourhood of $1$,
$F_{k_{0}}\left(y\right)\neq0$ and $F_{k}\left(y\right)=0$ for $0<k<k_{0}$ odd.
Assume that $x_{0}\left(y\right)$ and $F_{k_{0}}\left(y\right)$ are analytic at $y=1$, and that $x_{0}\left(y\right)$ satisfies
the \emph{variance condition}, $x_{0}''\left(1\right)x_{0}\left(1\right)+x_{0}'\left(1\right)x_{0}\left(1\right)-x_{0}'\left(1\right)^{2}\neq0$.

Let $X_{n}$ be the restriction of $\xi$ onto all objects in $\mathcal{F}$
of size $n$. Under these conditions, the distribution of $X_{n}$
is asymptotically Gaussian with mean \[
\mathbb{E}\left[X_{n}\right]\sim \mu n \ \ \ \mathrm{with}\ \ \mu=-\frac{x_{0}'\left(1\right)}{x_{0}\left(1\right)}%n+\frac{F_{k_{0}}'\left(1\right)}{F_{k_{0}}\left(1\right)}+O\left(\frac{1}{n}\right)
\]
 and variance\[
\mathbb{V}\left[X_{n}\right]\sim \sigma^2n\ \ \ \mathrm{with}\ \ \sigma^2=-\frac{x_{0}''\left(1\right)}{x_{0}\left(1\right)}-\frac{x_{0}'\left(1\right)}{x_{0}\left(1\right)}+\left(\frac{x_{0}'\left(1\right)}{x_{0}\left(1\right)}\right)^{2}%n+\left(\frac{F_{k_{0}}''\left(1\right)}{F_{k_{0}}\left(1\right)}+\frac{F_{k_{0}}'\left(1\right)}{F_{k_{0}}\left(1\right)}-\left(\frac{F_{k_{0}}'\left(1\right)}{F_{k_{0}}\left(1\right)}\right)^{2}\right)+O\left(\frac{1}{n}\right)
.\]

\end{thm}

\section{\label{sec:cycleindexsums}Exact enumeration of outerplanar graphs}
From now on we always 
consider outerplanar graphs as unlabeled objects, unless stated
otherwise.
In Section \ref{subsec:twoconn} we derive the cycle index sums for rooted and unrooted
two-connected outerplanar graphs.
Section \ref{subsec:con} shows how to decompose connected graphs
into two-connected components, which yields expressions of the cycle index sums
for rooted and unrooted connected outerplanar graphs.
In Section \ref{subsec:notcon} we use the simple fact
that an unrooted outerplanar graph is an (unordered) col\-lection of unrooted connected outerplanar graphs to obtain the cycle index sum for outerplanar graphs.
In Section \ref{subsec:bicon} we explain how to adapt the decomposition to enumerate bipartite outerplanar graphs.

\subsection{Enumeration of dissections (two-connected outerplanar graphs)}\label{subsec:twoconn}
A graph is \emph{two-connected} if at least two of its vertices have to be removed to disconnect it. 
A two-connected outerplanar graph with at least three vertices has
a unique Hamiltonian cycle \cite{LeydoldStadler} and can therefore
be embedded uniquely in the plane so that this Hamiltonian cycle
lies on the outer face. This unique embedding is thus a dissection of
a convex polygon. Hence the task of counting two-connected
outerplanar graphs coincides with the task of counting dissections
of a polygon. The generating functions for unlabeled (rooted
and unrooted) dissections were derived by Read \cite{Read}. In this
section, we extend Read's work to cycle index sums.

A dissection
can have two types of automorphisms: reflections
and rotations. We use the terminology of an \emph{oriented} dissection when an orientation
is imposed on the Hamiltonian cycle of the dissection. Since reflections
reverse orientations, they are excluded from the automorphism group
of an oriented dissection. By rooting at an edge on the
outer face, rotations are excluded from the automorphism group as well.
Thus, oriented dissections that are rooted at an edge on the outer
face are easy to count and are the starting point for the enumeration
of different types of dissections. Having the cycle index sum for
oriented and non-oriented dissections that are rooted at an edge on
the outer face, we can derive the cycle index sums for dissections
rooted at an edge not on the outer face, and for dissections rooted
at a face, by composition. Finally, we consider the dual of a dissection,
which is essentially a tree, and thereby count unrooted dissections. 
The dissimilarity characteristic
theorem for trees (Lemma \ref{lem:dissimilar_cycleindex})
can be applied to express the cycle index sum for
unrooted dissections by a combination of the cycle index sums for
several types of edge rooted and face rooted dissections.

We do not consider the graph consisting of a single vertex as a dissection.
However, it is con\-ven\-ient to include the single edge to the sets of edge rooted, vertex rooted and unrooted dissections.

\begin{terminology}[rooting, outer-edge, inner-edge, 
symmetry-edge, reflective] We say that a
dissection is \emph{edge} (respectively \emph{face}) \emph{rooted} if one edge (respectively face)
is distinguished from the others in the dissection. Note that the vertices of the
root edge (respectively face) might be interchanged
by the automorphisms of the dissection.
An edge on the
outer face is called \emph{outer-edge},
and \emph{inner-edge} otherwise.
 A \emph{symmetry-edge}
is an inner-edge such that there exists a nontrivial automorphism that fixes this
edge. It is clear that a dissection can have at most one symmetry-edge.
An edge rooted dissection is called \emph{reflective} if its automorphism
group contains a non-trivial reflection.
\end{terminology}

In the following we present the cycle index sums for several types of dissections and only sketch the proofs here. 
The corresponding notation is introduced
in Table \ref{cap:notationdissections}. 
We denote sets of graphs by calligraphic letters, 
ordinary generating functions by capital letters,
and counts by small letters. 
Details can be found in the thesis of the last author~\cite{Vigerske}.

\begin{table}[!ht]
\begin{center}\begin{tabular}{|ll|l|l|}
\hline 
\multicolumn{2}{|c|}{Notation}&
\multicolumn{2}{c|}{Type of dissection}\tabularnewline
\hline
\hline 
$\OEDc$&
$\OED\left(x\right)$&
oriented&
rooted at an outer-edge\tabularnewline
$\REDc$&
$\RED\left(x\right)$&
reflective&
rooted at an outer-edge\tabularnewline
$\IDc$&
$\ID\left(x\right)$&
not oriented&
rooted at an inner-edge\tabularnewline
$\SDc$&
$\SD\left(x\right)$&
not oriented&
rooted at a symmetry-edge\tabularnewline
\hline 
$\OFDc$&
$\OFD\left(x\right)$&
oriented&
rooted at an inner-face\tabularnewline
$\FDc$&
$\FD\left(x\right)$&
not oriented&
rooted at an inner-face\tabularnewline
$\SFDc$&
$\SFD\left(x\right)$&
not oriented&
rooted at an inner-face that is incident to a symmetry-edge\tabularnewline
\hline
$\VDc$&
$\VD\left(x\right)$&
not oriented&
rooted at a vertex\tabularnewline
\hline
$\Dc$&
$\D\left(x\right)$&
not oriented&
unrooted\tabularnewline
\hline
\end{tabular}\end{center}
\caption{\label{cap:notationdissections}Notation for several types of dissections}
\end{table}

We first recall the well-known results on oriented outer-edge rooted dissections (\cite{Read} and \cite[A001003]{Sloane}). A non-oriented
outer-edge rooted dissection is counted twice by the cycle index sum for oriented outer-edge rooted dissections when the dissection is not invariant under a reflection that fixes
the root-edge. Therefore, to derive counting formulas for non-oriented
dissections later, we also need the cycle index sum for reflective
outer-edge rooted dissections. Automorphisms of such structures can be
divided into two classes: Those that fix the vertices of the
root-edge and those that interchange the vertices of the root edge.

\begin{figure}[!ht]
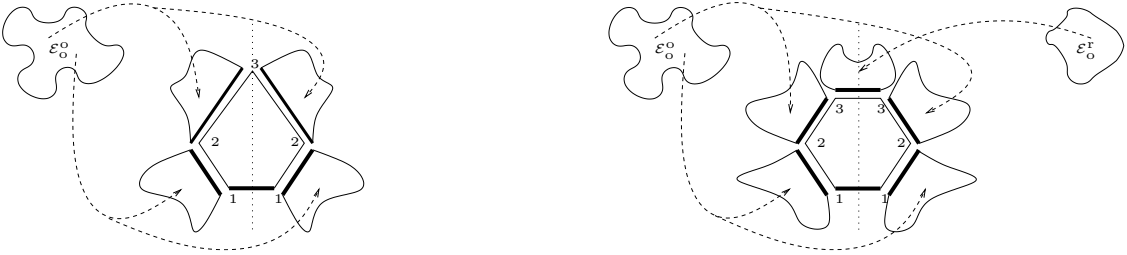

\begin{center}\input{red1_paper.pstex_t}\hfill{}\input{red2_paper.pstex_t}\end{center}

\caption{\label{cap:RED}Construction of a reflective outer-edge rooted dissection
for an odd and an even number of edges in the rooted polygon}
\end{figure}

\begin{lem}
[outer-edge rooted dissections]\label{lem:OED_RED}
The cycle index sum for oriented outer-edge rooted dissections is given by\begin{equation}
Z\left(\OEDc\right)=\frac{s_{1}}{4}\left(s_{1}+1-\sqrt{s_{1}^{2}-6s_{1}+1}\right).\label{eq:cycleindexOED}\end{equation}
 The cycle index sum for reflective outer-edge rooted dissections
is $Z\left(\REDc\right)=\frac{1}{2}\left(Z^{+}\left(\REDc\right)+Z^{-}\left(\REDc\right)\right),$
where $Z^{+}\left(\REDc\right)$ counts the mappings that fix the vertices of the root-edge (i.e., the identity mappings), and $Z^{-}\left(\REDc\right)$
counts the mappings that interchange the vertices of the root-edge
(i.e, the reflection mappings),\begin{eqnarray}
Z^{+}\left(\REDc\right) & = & \frac{1+s_{1}-3s_{1}^{2}+s_{1}^{3}-\left(1+s_{1}\right)\sqrt{s_{1}^{4}-6s_{1}^{2}+1}}{4s_{1}},\label{eq:cycleindexREDid}\\
Z^{-}\left(\REDc\right) & = & \frac{s_{2}^{2}-3s_{1}s_{2}+s_{2}+s_{1}-\left(s_{2}+s_{1}\right)\sqrt{s_{2}^{2}-6s_{2}+1}}{4s_{2}}.\label{eq:cycleindexREDr}\end{eqnarray}

\end{lem}
\begin{proof}
As pointed out before, oriented outer-edge rooted dissections have
a trivial automorphism group. Therefore, their numbers are closely related
to the bracketing numbers \cite[A001003]{Sloane}. To construct
such a dissection, we replace the edges of an edge rooted polygon other than the root-edge by oriented outer-edge
rooted dissections, see Figure \ref{cap:OED}. The cycle index
sum (\ref{eq:cycleindexOED}) is derived by an
application of the composition
formula (\ref{eq:composition}).

\begin{figure}[!ht]
\begin{center}\input{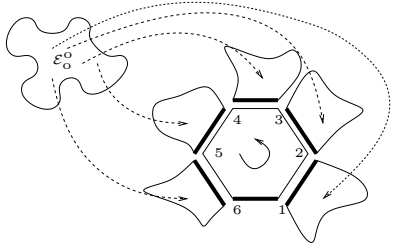}\end{center}

\caption{\label{cap:OED}Construction of an oriented outer-edge rooted dissection}
\end{figure}

Similarly to the construction of oriented outer-edge rooted dissections
a reflective outer-edge rooted dissection can be constructed by plugging
oriented outer-edge rooted dissections into the edges of an edge-rooted
polygon other than the root-edge. However, to obtain a reflective
dissection, edges of the polygon with the same distance to the root-edge
have to receive the same outer-edge rooted dissection.
If the polygon has an even number of edges,
the outer-edge rooted dissection that is plugged into the edge opposite to the root-edge has to be reflective itself, see also
Figure \ref{cap:RED}. 
 Corresponding formulas for the cycle index
sums are derived by an application of the composition formula (\ref{eq:composition}),\[
Z^{+}\left(\REDc\right)=\frac{s_{1}^{4}+\left(s_{1}-s_{1}^{2}\right)Z\left(\OEDc;s_{1}^{2}\right)}{s_{1}^{2}-2Z\left(\OEDc;s_{1}^{2}\right)},\qquad Z^{-}\left(\REDc\right)=\frac{s_{2}^{2}+\left(s_{1}-s_{2}\right)Z\left(\OEDc;s_{2}\right)}{s_{2}-2Z\left(\OEDc;s_{2}\right)}.\]
 Substituting formula (\ref{eq:cycleindexOED}) for $Z\left(\OEDc\right)$
yields formulas (\ref{eq:cycleindexREDid}) and (\ref{eq:cycleindexREDr}).
\end{proof}

\begin{lem}
[inner-edge and symmetry-edge rooted dissections]\label{lem:ID_SD}The
cycle index sums for inner-edge rooted dissections and symmetry-edge
rooted dissections, respectively, are given by\begin{eqnarray*}
Z\left(\IDc\right) & = & \frac{\left(3s_{1}-1+\sqrt{s_{1}^{2}-6s_{1}+1}\right)^{2}}{64}+\frac{\left(s_{1}+s_{2}\right)^{2}}{16s_{2}}\left(\frac{1+s_{2}-\sqrt{s_{2}^{2}-6s_{2}+1}}{1-s_{2}+\sqrt{s_{2}^{2}-6s_{2}+1}}\right)^{2}\\
 &  & -\frac{s_{1}^{2}+s_{2}}{16s_{2}}\left(3s_{2}-1+\sqrt{s_{2}^{2}-6s_{2}+1}\right),\\
Z\left(\SDc\right) & = & \frac{s_{1}^{6}-2s_{1}^{6}s_{2}^{2}+s_{1}^{4}s_{2}^{2}-s_{1}^{2}s_{2}^{3}-s_{2}^{3}}{16s_{1}^{4}s_{2}^{3}}+\frac{3}{8}\left(1-s_{2}-s_{1}^{2}\right)+\frac{1+s_{1}^{2}}{16s_{1}^{4}}\sqrt{s_{1}^{8}-6s_{1}^{4}+1}\\
&  & -\frac{1}{8}\sqrt{s_{1}^{4}-6s_{1}^{2}+1} -\frac{1}{16s_{2}}\left(1+\frac{s_{1}^{2}}{s_{2}^{2}}\right)\sqrt{s_{2}^{4}-6s_{2}^{2}+1}-\frac{1}{16}\left(1+\frac{s_{1}^{2}}{s_{2}}\right)\sqrt{s_{2}^{2}-6s_{2}+1}.\end{eqnarray*}

\end{lem}
\begin{proof}
An inner-edge rooted dissection can be constructed by joining two
outer-edge rooted dissections at their root-edge. We consider
the possible reflections and rotations. There are four kinds of
transformations of the plane that map the root-edge onto itself (also see
Figure \ref{cap:ID}): 1.~the identity mapping; 2.~the reflection
at the root-edge; 3. the half-turn around the root-edge; 4.~the reflection at the perpendicular bisector of the root-edge.
\begin{figure}[!ht]
\begin{center}\includegraphics{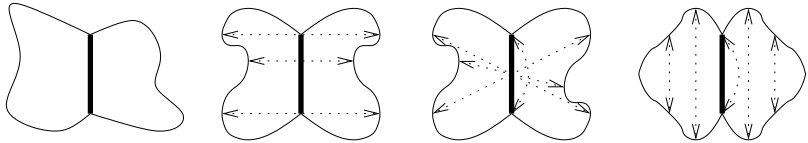}\end{center}
\caption{\label{cap:ID}The four kinds of transformations that fix the root
edge of an inner-edge rooted dissection.}
\end{figure}

Those inner-edge rooted dissections that are invariant under the identity
mapping are counted by joining two independently chosen non-empty
oriented outer-edge rooted dissections at their root-edge. Those that
are invariant under the second or third kind of permutation are constructed
by joining a non-empty oriented outer-edge rooted dissection with
a copy of itself on its root-edge, where the vertices of the root-edge
remain fixed (type 2) or are interchanged (type 3). Finally, inner-edge
rooted dissections invariant under the fourth kind of mapping are
composed out of two independently
chosen non-empty reflective outer-edge rooted dissections, joined
at the root-edge. Hence,
we get\[
Z\left(\IDc\right)=\frac{1}{4}\left(\left(\frac{Z\left(\OEDc\right)}{s_{1}^2}-1\right)^{2}s_1^2+\left(\frac{Z\left(\OEDc;s_{2}\right)}{s_{2}^{2}}-1\right)\left(s_{1}^{2}+s_{2}\right)+\left(\frac{Z^{-}\left(\REDc\right)}{s_{2}}-1\right)^{2}s_{2}\right),\]
 which together with Lemma \ref{lem:OED_RED} yields the statement.

Similarly, a non-empty dissection rooted at a symmetry-edge can be
constructed by joining an oriented outer-edge rooted dissection
with a copy of itself at the root-edge. We have to further
distinguish dissections that are invariant under reflection
and dissections that are not. A similar discussion as for inner-edge
rooted dissections leads to
%\begin{eqnarray*}
%Z\left(\SDc\right) & = & \frac{1}{4}\left(\frac{2}{s_{1}^{2}}\left(Z\left(\OEDc;s_{1}^{2}\right)-Z^{+}\left(\REDc;s_{1}^{2},s_{2}^{2}\right)\right)+\frac{s_{1}^{2}+s_{2}}{s_{2}^{2}}\left(Z\left(\OEDc;s_{2}\right)-Z^{+}\left(\REDc;s_{2},s_{4}\right)\right)\right)\\
%&  & +\frac{1}{4}\left(\frac{1}{s_{1}^{2}}Z^{+}\left(\REDc;s_{1}^{2},s_{2}^{2}\right)-s_{1}^{2}+\frac{s_{1}^{2}}{s_{2}^{2}}Z^{+}\left(\REDc;s_{2},s_{4}\right)-s_{1}^{2}\right)\\
%&  & +\frac{1}{4}\left(\frac{1}{s_{2}}Z^{-}\left(\REDc;s_{1}^{2},s_{2}^{2}\right)-s_{2}+\frac{1}{s_{2}}Z^{-}\left(\REDc;s_{2},s_{4}\right)-s_{2}\right),\end{eqnarray*}
%which can be simplified to
\[
Z\left(\SDc\right)=\frac{Z\left(\OEDc;s_{1}^{2}\right)}{2s_{1}^{2}}+\frac{s_{1}^{2}+s_{2}}{4s_{2}^{2}}Z\left(\OEDc;s_{2}\right)-\frac{Z^{+}\left(\REDc;s_{1}^{2},s_{2}^{2}\right)}{4s_{1}^{2}}+\frac{Z^{-}\left(\REDc;s_{1}^{2},s_{2}^{2}\right)}{4s_{2}}-\frac{s_{1}^{2}+s_{2}}{2}.\]
We then obtain the statement by applying Lemma \ref{lem:OED_RED}.
\end{proof}

We proceed with dissections that are rooted at a face. In this case,
the automorphism groups might include cyclic permutations of order
greater than two.

\begin{figure}[!ht]
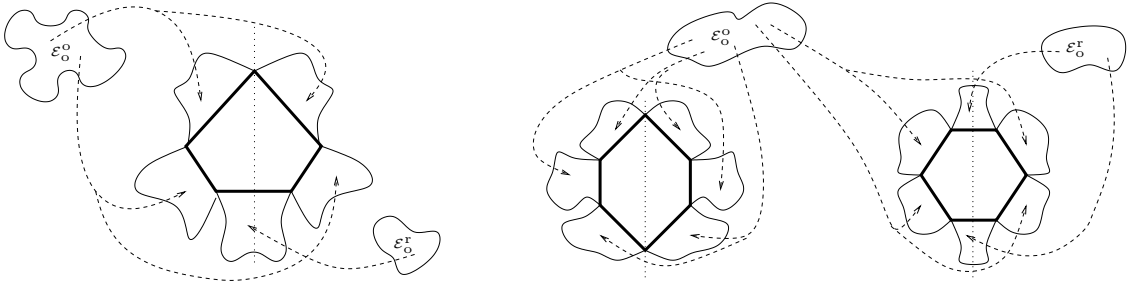

\begin{center}\input{fd1_paper.pstex_t}\hfill{}\input{fd2_paper.pstex_t}\end{center}
\caption{\label{cap:FD}The extra terms $s_{1}s_{2}^{m}$ ($k=2m+1$, left)
and $s_{2}^{m+1}$ and $s_{1}^{2}s_{2}^{m}$ ($k=2m+2$, right) in
$Z\left(D_{k}\right)$ for $m=2$.}
\end{figure}

\begin{lem}[face rooted dissections]\label{lem:FD}
The cycle index sum for face
rooted dissections is given by\begin{multline*}
Z\left(\FDc\right)=-\frac{1}{2}\sum_{d\geq1}\frac{\varphi\left(d\right)}{d}\log\left(\frac{3}{4}-\frac{1}{4}s_{d}+\frac{1}{4}\sqrt{s_{d}^{2}-6s_{d}+1}\right)+\frac{s_{1}+5}{32}\sqrt{s_{1}^{2}-6s_{1}+1}\\
-\frac{s_{1}^{2}+2s_{1}+2s_{2}+7}{32}+\frac{s_{1}^{2}+6s_{1}s_{2}+8s_{1}^{2}s_{2}+3s_{2}^{2}+4s_{1}s_{2}^{2}+5s_{1}^{2}s_{2}^{2}+2s_{2}^{3}-2s_{1}s_{2}^{3}-2s_{1}^{2}s_{2}^{3}-s_{2}^{4}}{16s_{2}\left(1-s_{2}+\sqrt{s_{2}^{2}-6s_{2}+1}\right)^{2}}\\
+\left(\frac{1}{16}-\frac{s_{1}^{2}+2s_{1}s_{2}+3s_{1}^{2}s_{2}+s_{2}^{2}-2s_{1}s_{2}^{2}-2s_{1}^{2}s_{2}^{2}-s_{2}^{3}}{8s_{2}\left(1-s_{2}+\sqrt{s_{2}^{2}-6s_{2}+1}\right)^{2}}\right)\sqrt{s_{2}^{2}-6s_{2}+1}.\end{multline*}
\end{lem}
\begin{proof}
A face rooted dissection can be constructed by plugging oriented outer-edge
rooted dissections into the edges of a polygon. In the oriented case (see Figure \ref{cap:OFD}), we only have to consider cyclic permutations
of the polygon. Let $C_{k}$ be the cyclic group,
generated by the permutation $\left(1\,2\,3\,\cdots\, k\right)$.
The composition formula (\ref{eq:composition}) gives
the cycle index sum for oriented face rooted dissections,
\begin{equation}
Z\left(\OFDc\right)=\sum_{k\geq3}Z\left(C_{k}\right)\left[Z\left(\OEDc\right)/s_{1}\right].\label{eq:cycleindexsumOFDintermediate}\end{equation}
\begin{figure}[!ht]
\begin{center}\input{ofd_paper.pstex_t}\end{center}
\caption{\label{cap:OFD}Construction of a face rooted dissection counted
by $s_{3}^{2}\left[Z\left(\OEDc\right)/s_{1}\right]$}
\end{figure}
$Z\left(C_{k}\right)$ can be expressed with the Euler-$\varphi$-function,
$Z\left(C_{k}\right)=\frac{1}{k}\sum_{d|k}\varphi\left(d\right)s_{d}^{k/d},$ which leads to
\[Z\left(\OFDc\right)=-\left\{ \sum_{d\geq1}\frac{\varphi\left(d\right)}{d}\log\left(1-\frac{Z\left(\OEDc;s_{d}\right)}{s_{d}}\right)\right\} -\frac{Z\left(\OEDc\right)}{s_{1}}-\frac{1}{2}\left(\left(\frac{Z\left(\OEDc\right)}{s_{1}}\right)^{2}+\frac{Z\left(\OEDc;s_{2}\right)}{s_{2}}\right).\]

In the non-oriented case, we have to take care of
additional reflections. Therefore, the cyclic group $C_{k}$ has to
be replaced by the dihedral group $D_{k}$, which is generated by
the cycle $\left(1\,2\,3\,\cdots\, k\right)$ and the reflection $\left(1\: k\right)\left(2\:(k-1)\right)\left(3\:(k-2)\right)\cdots$, and has the cycle index\[
Z\left(D_{k}\right)=\frac{1}{2}Z\left(C_{k}\right)+\begin{cases}
\frac{1}{2}s_{1}s_{2}^{m}, & k\textrm{ odd},\ k=2m+1, \\
\frac{1}{4}s_{2}^{m+1}+\frac{1}{4}s_{1}^{2}s_{2}^{m}, & k\textrm{
even},\ k=2m+2.\end{cases}\]
The objects with cyclic automorphisms are counted by $Z\left(C_{k}\right)$.
For the reflections, the outer-edge
rooted dissections attached to corresponding pairs of edges must be
the same, while the outer-edge rooted dissection that is mapped to
itself must be reflective. We also have to distinguish
between polygons of odd
and even size. If we identify the corresponding terms
in $Z\left(D_{k}\right)$
with the correct cycle index sums (see Figure \ref{cap:FD}) we get \[
Z\left(\FDc\right)=\frac{1}{2}Z\left(\OFDc\right)+\frac{1}{2s_{2}}\left(s_{1}Z^{-}\left(\REDc\right)+\frac{s_{1}^{2}}{2s_{2}}Z\left(\OEDc;s_{2}\right)+\frac{1}{2}Z^{-}\left(\REDc\right)^{2}\right)\:\frac{Z\left(\OEDc;s_{2}\right)}{s_{2}-Z\left(\OEDc;s_{2}\right)},\]
which together with Lemma \ref{lem:OED_RED} yields the statement.
\end{proof}

\begin{lem}[face rooted dissections with root-face incident to a symmetry-edge]
\label{lem:SFD}The cycle index sum for face rooted dissections where
the root-face is incident to a symmetry-edge is given by\begin{multline*}
Z\left(\SFDc\right)=\frac{s_{1}^{6}\left(1-3s_{2}^{2}-3s_{2}^{3}\right)+s_{1}^{4}\left(s_{2}^{2}+5s_{2}^{3}-3s_{2}^{4}\right)-s_{2}^{3}-s_{1}^{2}s_{2}^{3}}{8s_{1}^{4}s_{2}^{3}}\\
-\frac{1}{4}\sqrt{s_{1}^{4}-6s_{1}^{2}+1}+\frac{s_{2}^{3}+s_{1}^{2}s_{2}^{3}}{8s_{1}^{4}s_{2}^{3}}\sqrt{s_{1}^{8}-6s_{1}^{4}+1}-\frac{s_{1}^{2}+s_{2}^{2}}{8s_{2}^{3}}\sqrt{s_{2}^{4}-6s_{2}^{2}+1}.\end{multline*}

\end{lem}
\begin{proof}
As in the case of symmetry-edge rooted dissections 
we join two (identical) non-empty oriented
outer-edge rooted dissections at their root-edge and choose one of
the faces incident to the root-edge to be the root face. In contrast
to Lemma \ref{lem:ID_SD} we do not have to care about the permutations
of the second and third kind, and obtain

\[
Z\left(\SFDc\right)=\frac{1}{s_{1}^{2}}Z\left(\OEDc;s_{1}^{2}\right)-\frac{1}{2s_{1}^{2}}Z^{+}\left(\REDc;s_{1}^{2},s_{2}^{2}\right)+\frac{1}{2s_{2}}Z^{-}\left(\REDc;s_{1}^{2},s_{2}^{2}\right)-\frac{1}{2}\left(s_{1}^{2}+s_{2}\right),\]
which together with Lemma \ref{lem:OED_RED} yields the statement.
\end{proof}
We have now found all cycle index sums that are needed to compute
the cycle index sum for unrooted dissections.

\begin{thm}
[dissections]\label{thm:Dcycleindex}The cycle index sum for dissections
is given by\begin{align*}
Z\left(\Dc\right) = & -\frac{1}{2}\sum_{d\geq1}\frac{\varphi\left(d\right)}{d}\log\left(\frac{3}{4}-\frac{1}{4}s_{d}+\frac{1}{4}\sqrt{s_{d}^{2}-6s_{d}+1}\right)+\frac{s_{2}+s_{1}^{2}-4s_{1}-2}{16}\\
 & +\frac{s_{1}^{2}-3s_{1}^{2}s_{2}+2s_{1}s_{2}}{16s_{2}^{2}}+\frac{3-s_{1}}{16}\sqrt{s_{1}^{2}-6s_{1}+1}-\frac{1}{16}\left(1+\frac{s_{1}^{2}}{s_{2}^{2}}+\frac{2s_{1}}{s_{2}}\right)\sqrt{s_{2}^{2}-6s_{2}+1}.\end{align*}

\end{thm}
\begin{proof}
Since the dual graph of a dissection that has at least one face is
a tree, we can apply the dissimilarity characteristic theorem for
trees in its cycle index version (Lemma \ref{lem:dissimilar_cycleindex})
to derive the cycle index sum for unlabeled dissections from the ones
for face rooted dissections (corresponding to vertex rooted trees),
inner-edge rooted dissections (corresponding to edge rooted trees),
and symmetry-edge rooted dissections (corresponding to symmetry-edge
rooted trees). A vertex of the tree that is not incident to a symmetry-edge
corresponds to a face of the dissection that is not incident to a
symmetry-edge. Lemma~\ref{lem:dissimilar_cycleindex} with $\FDc\setminus\SFDc$
instead of $\mathcal{V}_{1}$, $\IDc$ instead of $\mathcal{E}$, and $\SDc$ instead of $\mathcal{S}$ yields\[
Z\left(\Dc\right)=\frac{1}{2}\left(s_{1}^{2}+s_{2}\right)+Z\left(\FDc\right)-Z\left(\SFDc\right)-Z\left(\IDc\right)+2Z\left(\SDc\right).\]
The additional term $\frac{1}{2}\left(s_{1}^{2}+s_{2}\right)$ counts
the dissection that consists of one edge only. We apply Lemmas \ref{lem:ID_SD}, \ref{lem:FD}, and~\ref{lem:SFD} to obtain the result.
\end{proof}

Replacing $s_1$ by $x$, $s_2$ by $x^2$, $\ldots$ we obtain the generating function $D(x)$ of
dissections, which was already found by Read:
\begin{align}
\D\left(x\right) = & -\frac{1}{2}\sum_{d\geq1}\frac{\varphi\left(d\right)}{d}\log\left(\frac{1}{4}\left(3-x^{d}+\sqrt{x^{2d}-6x^{d}+1}\right)\right)\label{eq:Dogf}\\
 &
+\frac{x^{2}}{8}-\frac{1}{4}x-\frac{5}{16}+\frac{1}{8x}+\frac{1}{16x^{2}}+\frac{3-x}{16}\sqrt{x^{2}-6x+1}-\frac{1+2x+x^{2}}{16x^{2}}\sqrt{x^{4}-6x^{2}+1}.\nonumber
\end{align}
The coefficients of $D(x)$, counting unlabeled dissections, can be extracted in polynomial time, $D(x)=x^2+x^3+2x^4+3x^5+9x^6+20x^7+75x^8+262x^9+\ldots$, matching the values computed by Read, see \cite[A001004]{Sloane}.

Finally, the cycle index sum for vertex rooted dissections, which we will need in Section \ref{subsec:con}, can be
derived by using Formula (\ref{eq:rooting}):
$Z\left(\VDc\right)=s_{1}\frac{\partial}{\partial s_{1}}Z\left(\Dc\right)$.

\begin{cor}
[vertex rooted dissections]\label{cor:VD}The cycle index sum for
vertex rooted dissections is given by\begin{equation}\label{eq:VDcycleindex}
%Z\left(\VDc\right)=\frac{s_{1}}{8s_{2}^{2}}\left(s_{2}+s_{1}-3s_{1}s_{2}+s_{1}s_{2}^{2}-2s_{2}^{2}-\left(s_{1}+s_{2}\right)\sqrt{s_{2}^{2}-6s_{2}+1}-s_{2}^{2}\sqrt{s_{1}^{2}-6s_{1}+1}\right).
Z\left(\VDc;s_1,s_2\right)=\frac{s_{1}}{8}\left(1+s_1-\sqrt{s_1^2-6s_1+1}\right)
+\frac{s_1}{8s_2^2}(s_1+s_2)\left(1-3s_2-\sqrt{s_2^2-6s_2+1}\right).
\end{equation}
\end{cor}

\subsection{Enumeration of connected outerplanar graphs}\label{subsec:con}
We denote the set of unrooted connected outerplanar graphs by $\Cc$, 
and the set of vertex rooted connected outerplanar graphs by $\VCc$.
All rooted graphs considered in this section are rooted at a vertex. Again, ordinary generating functions are denoted
by capital letters and coefficients by small letters. Thus, $\VC\left(x\right)=\sum_{n}\vc_{n}x^{n}$
and $\C\left(x\right)=\sum_{n}\con_{n}x^{n}$.

The cycle index sum for rooted connected outerplanar graphs is derived
by decomposing the graphs into rooted two-connected outerplanar graphs,
i.e., vertex rooted dissections.

\begin{lem}
[rooted connected outerplanar graphs]\label{lem:VCcycleindex}The
cycle index sum for vertex rooted connected outerplanar graphs is
implicitly determined by the equation\begin{equation}
Z(\VCc)=s_{1}\exp\left(\sum_{k\geq1}\frac{Z\left(\VDc;Z\left(\VCc;s_{k},s_{2k},\ldots\right),Z\left(\VCc;s_{2k},s_{4k},\ldots\right)\right)}{k\, Z\left(\VCc;s_{k},s_{2k},\ldots\right)}\right).\label{eq:VCcycleindex}\end{equation}

\end{lem}
\begin{proof}
Graphs in $\VCc$ rooted at a vertex that is not a cut-vertex
can be constructed by taking a rooted dissection and attaching a rooted
connected outerplanar graph at each vertex of the dissection other
than the root vertex. By the composition formula (\ref{eq:composition})
we obtain that
\begin{equation}
s_{1}\left(\frac{Z\left(\VDc\right)}{s_{1}}\right)\left[Z(\VCc)\right]\label{eq:noncut-vertexrootedblock}
\end{equation}
is the cycle index sum for connected outerplanar graphs rooted at a non-cut-vertex.
The division (resp. multiplication) by $s_{1}$ is due to the removal (resp. addition) of the root vertex before (resp. after)
application of Formula (\ref{eq:composition}).

The cycle index sum for rooted connected outerplanar graphs
where the root vertex is incident to exactly $n$ blocks, $n\geq2$,
can be obtained by another application of the composition theorem.
We join $n$ connected outerplanar graphs that are rooted at a vertex
other than a cut-vertex at their root vertex. Application of the composition
formula (\ref{eq:composition}) with the symmetric group $S_{n}$
and Formula (\ref{eq:noncut-vertexrootedblock}) (divided by $s_{1}$)
for the cycle index sum for non-cut-vertex rooted connected outerplanar graphs (excluding the root) yields\[
s_{1}Z\left(S_{n}\right)\left[\left(\frac{Z(\VDc)}{s_{1}}\right)\left[Z(\VCc)\right]\right].\]
Summing over $n\geq 0$, we get
($Z\left(S_{0}\right):=1$)
\[
Z(\VCc) =
s_{1}\sum_{n\geq0}Z\left(S_{n}\right)\left[\left(\frac{Z\left(\VDc\right)}{s_{1}}\right)\left[Z(\VCc)\right]\right].
\]
With the well-known formula $\sum_{n\geq0}Z\left(S_{n}\right)=\exp\left(\sum_{k\geq1}\frac{1}{k}s_{k}\right)$, the statement follows.
\end{proof}
\begin{thm}
[connected outerplanar graphs]\label{thm:Ccycleindex}The cycle index
sum for connected outerplanar graphs is given by\begin{equation}
Z\left(\Cc\right)=Z(\VCc)+Z\left(\Dc;Z(\VCc)\right)-Z\left(\VDc;Z(\VCc)\right).\label{eq:Ccycleindex}\end{equation}

\end{thm}
\begin{proof}
To derive the cycle index sum for unrooted connected outerplanar graphs,
one can use Formula (\ref{eq:unrooting}). Hence,\begin{equation}
Z\left(\Cc\right)=\int_{0}^{s_{1}}\frac{1}{s_{1}}Z(\VCc)ds_{1}+Z\left(\Cc\right)\left|_{s_{1}=0}\right..\label{eq:Cbyintegration}\end{equation}
The term $Z\left(\Cc\right)\left|_{s_{1}=0}\right.$ can be further
replaced by $Z\left(\Dc\right)\left|_{s_{1}=0}\right.[Z(\VCc)]$ 
because each fixed-point free permutation in a connected
graph $G$ has a unique block whose vertices are setwise fixed by 
the automorphisms of $G$ \cite[page 190]{HararyPalmer}. 
Using the special structure (\ref{eq:VCcycleindex}) of $Z(\VCc)$,
a closed solution for the integral in (\ref{eq:Cbyintegration}) can be found \cite{Vigerske}.
We put these facts together and obtain (\ref{eq:Ccycleindex}).
\end{proof}
Replacing $s_i$ by $x^i$ in $Z(\VCc)$, we obtain that the generating function $\hat{C}(x)$ counting vertex
rooted connected outerplanar graphs satisfies

\begin{equation}
\VC\left(x\right) = x\exp\left(\sum_{k\geq1}\frac{Z\left(\VDc;\VC\left(x^{k}\right)\right)}{k\,\VC\left(x^{k}\right)}\right),\label{eq:VCogf}
\end{equation}
from which the coefficients $\hat{C}_n$ counting vertex rooted
connected outerplanar graphs can be extracted in polynomial time:
$\VC\left(x\right)=x+x^2+3x^3+10x^4+40x^5+181x^6+918x^7+\ldots$, see
\cite{Tomii, Vigerske} for more entries. The numbers in \cite{Tomii}
verify the correctness of our result and were computed by the polynomial algorithm proposed in~\cite{BodirskyKang}.

In addition, it follows from~(\ref{eq:Ccycleindex}) that the generating function $C(x)$
counting connected outerplanar graphs satisfies:
\begin{equation}
\C\left(x\right)=\VC\left(x\right)+Z\left(\Dc;\VC\left(x\right)\right)-Z\left(\VDc;\VC\left(x\right)\right),\label{eq:Cogf}
\end{equation}
from which the coefficients $c_n$ counting connected outerplanar
graphs can be extracted in polynomial time:
$C(x)=x+x^2+2x^3+5x^4+13x^5+46x^6+172x^7+\ldots$, see \cite[A111563]{Sloane} for more entries.

\subsection{Enumeration of outerplanar graphs}\label{subsec:notcon}

We denote the set of outerplanar graphs by $\Outc$, its ordinary
generating function by $\Out\left(x\right)$ and the number of outerplanar
graphs with $n$ vertices by $\out_{n}$. As an
outerplanar graph is a collection of connected outerplanar graphs,
it is now easy to obtain the cycle index sum for outerplanar graphs.
An application of the composition formula (\ref{eq:composition})
with the symmetric group $S_{l}$ and object set $\Cc$ yields that $Z\left(S_{l}\right)\left[Z\left(\Cc\right)\right]$
is the cycle index sum for outerplanar graphs with $l$ connected components.
Thus, by summation over all $l\geq 0$ (we include here also the empty graph into $\Outc$ for convenience), we obtain the following theorem.

\begin{thm}
[outerplanar graphs]\label{thm:Outcycleindex}The cycle index sum
for outerplanar graphs is given by\[
Z\left(\Outc\right)=\exp\left(\sum_{k\geq1}\frac{1}{k}Z\left(\Cc;s_{k},s_{2k},\ldots\right)\right).\]

\end{thm}

Hence the generating functions $G(x)$ and $C(x)$ of outerplanar and
connected outerplanar graphs are related by
\begin{equation}
\Out\left(x\right)=\exp\left(\sum_{k\geq1}\frac{1}{k}\C\left(x^{k}\right)\right).\label{eq:Outogf}
\end{equation}
From this, we can extract in polynomial time the coefficients counting
outerplanar graphs,
$G(x)=1+x+2x^2+4x^3+10x^4+25x^5+80x^6+277x^7+\ldots$, see \cite[A111564]{Sloane} for more entries.

\subsection{Enumeration of bipartite outerplanar graphs}\label{subsec:bicon}
To study the chromatic number of a typical outerplanar graph
we enumerate bipartite outerplanar
graphs. Observe that an outerplanar graph is bipartite if and only
if all of its blocks are bipartite. As discussed
in Section \ref{sec:cycleindexsums},
blocks of an outerplanar graph are dissections, and it is clear that
a dissection is bipartite when all of its inner faces have an even
number of vertices. The decomposition of dissections exposed in
Section~\ref{subsec:twoconn} can be adapted to dissections where 
all faces have even degree. 
Once the cycle index sum for bipartite
dissections is obtained, the computation of the cycle index sums for bipartite connected
outerplanar graphs, and then of bipartite outerplanar graphs works in the same
way as for the general case, see~\cite{Vigerske} for details.
From that the coefficients of the series $G_{\textrm \it{b}}(x)$ counting bipartite
outerplanar graphs can be extracted in polynomial time: $G_{\textrm \it{b}}(x)=1+x+x^2+x^3+7x^4+12x^5+29x^6+61x^7+\ldots$, see the sequences A111757, A111758, and A111759 of
\cite{Sloane} for the coefficients of two-connected,
connected, and general bipartite outerplanar graphs.

\section{\label{sec:asymptotics}Asymptotic enumeration of unlabeled
outerplanar graphs}
To determine the asymptotic number of two-connected, connected, and
general outerplanar graphs, we use singularity analysis
as introduced in Section \ref{sub:singularityanalysis}. 
To compute the growth constants and subexponential
factors we expand the generating functions for outerplanar graphs 
around their dominant singularities.
For unlabeled two-connected
outerplanar graphs we present an analytic expression of the growth
constant. For the connected and the general case 
we give numerical approximations of the growth constants 
in Section~\ref{sec:approxgrowth}.

\subsection{Asymptotic estimates}\label{sec:estimates}
We now prove the first part of Theorem \ref{thm:asymptoticouter} on the asymptotic number of dissections. 

\begin{thm}[asymptotic number of unrooted dissections]\label{thm:twoconn}
The number $\dis_n$ of unlabeled two-connected outerplanar graphs on $n$ vertices has the asymptotic estimate
$\dis_{n} \sim \dis\, n^{-\frac{5}{2}}\delta^{-n}$
with growth rate $\delta^{-1}=3+2\sqrt{2}\approx5.82843$ and constant $\dis \approx0.00596026$.
\end{thm}

\begin{proof}
Let $\delta$ be the smallest root of $x^2-6x+1$, $\delta=3-2\sqrt{2}$. Equation (\ref{eq:Dogf}) implies that $D(x)$ can be
written as 
\[
\D\left(x\right)=-\frac{1}{2}\log\left(1-\frac{\sqrt{x^{2}-6x+1}}{x-3}\right)+\frac{3-x}{16}\sqrt{x^{2}-6x+1}+A\left(x\right),\]
 where $A\left(x\right)$ is analytic at $0$ with radius of convergence
$>\delta$. Since the logarithmic term is analytic for $\left|x\right|<\delta$,
we can expand it and collect ascending powers of $\sqrt{x^{2}-6x+1}$
in $\D\left(x\right)$. Thus,\[
\D\left(x\right)=\left(-\frac{1}{16\left(x-3\right)}+\frac{1}{6\left(x-3\right)^{3}}\right)\left(\sqrt{x^{2}-6x+1}\right)^{3}+\sum_{k\geq4}\frac{1}{2k}\left(\frac{\sqrt{x^{2}-6x+1}}{x-3}\right)^{k}+\tilde{A}\left(x\right),\]
 where $\tilde{A}\left(x\right)$ is again analytic at $0$ with radius
of convergence $>\delta$. Finally, using $\sqrt{x^{2}-6x+1}=\sqrt{1-x/\delta}\sqrt{1-\delta x}$
for $x\leq\delta$ and applying Lemma \ref{lem:basicscale}
we obtain\begin{align*}
\dis_{n} & =\left(-\frac{1}{16\left(\delta-3\right)}+\frac{1}{6\left(\delta-3\right)^{3}}\right)\left(\sqrt{1-\delta^{2}}\right)^{3}\frac{1}{\Gamma\left(-3/2\right)}n^{-5/2}\delta^{-n}\left(1+O\left(\frac{1}{n}\right)\right)\\
 & \sim\frac{\left(3\sqrt{2}-4\right)^{3/2}}{8\sqrt{2\pi}}n^{-5/2}\left(3+2\sqrt{2}\right)^{n}.\qedhere\end{align*}
\end{proof}

We now turn to the problem of asymptotic enumeration of connected
outerplanar graphs. First we have to establish the singular
development of the generating function for vertex rooted con\-nec\-ted outerplanar graphs $\VC\left(x\right)$.
Let $\rho$ be the radius of convergence of $\VC\left(x\right)$.
Observe that the coefficients $\vc_n$ are bounded from below by the number of unlabeled vertex rooted dissections $\vd_n$, 
which have exponential growth.
The coefficients are bounded from above 
by the number of embedded outerplanar graphs with a root edge,
which also have ex\-po\-nen\-tial growth 
(this follows from classical enumerative results 
on planar maps; see \cite{Tutte}).
Hence $\rho$ is in $(0,1)$.

To apply Theorem \ref{thm:positiveimplicitfunctions}
for rooted connected outerplanar graphs,
we consider the function\begin{equation}
H\left(x,y\right):=x\exp\left(\frac{Z\left(\VDc;y,\VC\left(x^{2}\right)\right)}{y}+\sum_{k\geq2}\frac{Z\left(\VDc;\VC\left(x^{k}\right),\VC\left(x^{2k}\right)\right)}{k\,\VC\left(x^{k}\right)}\right)-y.\label{eq:HforVC}
\end{equation}

Observe that Equation~(\ref{eq:VCogf}) implies that $H(x,\VC(x))=0$.
The difficulty in the application of the singular implicit functions theorem (Thm. \ref{thm:positiveimplicitfunctions})
is the verification of the requirements of this theorem.
Hence, to apply Theorem \ref{thm:positiveimplicitfunctions}, we have
to check that the do\-mi\-nant singularity of the generating functions for the connected components 
is determined by its im\-pli\-cit definition (like (\ref{eq:VCogf})) and not by a singularity of $H(x,y)$. This analysis is the main purpose of the next proposition. Observe that it can also be easily 
generalized to other classes of connected unlabeled graphs with known blocks.

\begin{lem}\label{lem:VCsingularity}
The generating function $\VC(x)$ satisfies the conditions of Theorem \ref{thm:positiveimplicitfunctions} with
the function $H\left(x,y\right)$ from Equation~(\ref{eq:HforVC}) and
$(r,s)=(\rho,\tau)$, where $\rho$ is the dominant singularity of $\VC(x)$ and $\tau:=\lim_{x\rightarrow \rho^-}\VC(x)$.

As a consequence, Theorem \ref{thm:positiveimplicitfunctions} ensures
that $\VC(x)$ has a singular expansion\begin{equation}
\VC\left(x\right)=\sum_{k\geq0}\VC_{k}X^{k},\quad\textrm{with}\quad X:=\sqrt{1-\frac{x}{\rho}},\quad\VC_{0}=\tau,\quad\VC_{1}=-\sqrt{\frac{2\rho\frac{\partial}{\partial x}H\left(\rho,\tau\right)}{\frac{\partial^{2}}{\partial y^{2}}H\left(\rho,\tau\right)}},\label{eq:VCasympt}
\end{equation}
with constants $\VC_{k}$, $k\geq2$, which can be computed from the
derivatives of $H\left(x,y\right)$ at $\left(\rho,\tau\right)$.

\end{lem}
\begin{proof}
The conditions $H(0,0)=0$ and $\frac{\partial}{\partial y}H(0,0)=-1$ can be verified easily. The positivity conditions on the coefficients of $H(x,y)$ follow from the positivity of the coefficients of $Z\left(\VDc\right)$.
The analyticity domain of $H(x,y)$ is determined
by the dominant singularities of $Z\left(\VDc\right)$; that is, $H\left(x,y\right)$
is analytic for $x$ and $y$ such that $\left|y\right|<\delta$ and $|x^l|<\rho$ and $|\VC\left(x^{l}\right)|<\delta$ for each
$l\geq2$. Since $\VC\left(x\right)$ is strictly increasing for positive
$x$ and since $\rho <1$, $|\VC\left(x^{l}\right)|\leq|\VC\left(x^{2}\right)|$ for all
$l\geq2$ and $|x|<\rho$. Therefore, $H\left(x,y\right)$ is analytic for $\left|x\right|<R:=\min(\sqrt{\rho},\sqrt{\VC^{-1}\left(\delta\right)})$
and $\left|y\right|<S:=\delta$.

We show next that $\rho < R$ and $\tau < S$.
\begin{enumerate}
\item We show $\tau\leq S=\delta$.
Let $\tilde H(x,y):=H(x,y)+y$. $\tilde H(x,y)$ satisfies $\tilde H(x,\VC(x))=\VC(x)$ and has the same domain of analyticity as $H(x,y)$.
Assume $\tau>\delta$. Then there exists $x_0<\rho$ such that $\VC\left(x_0\right)=\delta$.
Observe that, if $|x|< x_0$ then $|\VC(x^2)|\leq|\VC(x)|<\VC(x_0)=\delta$. Thus $(x,\VC(x))$ is in the
analyticity domain of $\tilde H(x,y)$, so that $\tilde H(x,\VC(x))=\VC(x)$.
By continuity we obtain $\tilde H(x_0,\VC(x_0))=\VC(x_0)$.
We have now the contradiction that
$\VC(x)$ is analytic at $x_0$ since $x_0<\rho$, whereas $\tilde H(x,\VC(x))$ is singular at $x_0$ because $\VC(x_0)=\delta$.

\item From 1 we know that $\tau\leq \delta$, i.e., $\rho\leq
\hat{C}^{-1}(\delta)$. Hence $R=\sqrt{\rho}>\rho$.

\item It remains to prove that $\tau<S$. Assume $\tau=\delta$.
Observe from (\ref{eq:VDcycleindex}) and (\ref{eq:VCogf}) that \[
\VC(x)=x \exp(\Psi(\VC(x))+A(x))
\]
where $\Psi(y)=1/8(1+y-\sqrt{1-6y+y^2})$ has a dominant singularity at $y=\delta$,
and where $A(x)$ is a generating function analytic for $|x|<\rho$ and having nonnegative coefficients.
(This follows from the fact that $2A(x)$ is the generating function for reflective vertex rooted dissections \cite{Vigerske}.)
Hence, for $0<x<\rho$,\[
\VC'(x)\geq \VC'(x)\Psi'(\VC(x))\VC(x),
\]
so that $\Psi'(\VC(x))\leq1/\VC(x)$. Thus, $\Psi'(\VC(x))$ is bounded when $x\rightarrow\rho^-$, which contradicts
the fact that $\lim_{y\rightarrow \delta^-}\Psi'(y)=+\infty$.
\end{enumerate}

Thus, $H(x,y)$ is analytic at $(\rho,\tau)$ and $H(\rho,\tau)=0$ is satisfied.
As pointed out before, the dominant singularity $\rho$ of $\VC(x)$ is determined either by a singularity in a component of
Equation~(\ref{eq:VCogf}), or by a non-uniqueness in the definition of $\VC(x)$ by Equation~(\ref{eq:VCogf}).
The relation $\tau<\delta$ excludes the first case, so that the singularity is caused by a non-uniqueness of the inversion.
Hence, the derivative of $H(x,y)$ with respect to $y$ has to vanish at $(x,y)=(\rho,\tau)$, since otherwise the implicit
function theorem ensures a (unique) analytic continuation of $\VC(x)$ at $x=\rho$.
Therefore, the equations from~(\ref{eq:posimplicitfunc}) are satisfied for $(r,s)=(\rho,\tau)$.

Furthermore, it is easily verified that
\begin{align*}
\left.\frac{\partial^2}{\partial y^2}H(x,y)\right|_{(x,y)=(\rho,\tau)} = &
\frac{1}{\tau}+\left.\frac{\partial^2}{\partial s_1^2}\frac{Z(\VDc;s_1,\VC(\rho^2))}{s_1}\right|_{s_1=\tau}
 = \frac{1}{\tau}+\frac{\tau}{\left(\tau^2-6\tau+1\right)^{3/2}}, \\
\left.\frac{\partial}{\partial x}H(x,y)\right|_{(x,y)=(\rho,\tau)} = & \tau\left(\frac{1}{\rho}+\left.\frac{\partial}{\partial x}
\frac{Z\left(\VDc;\tau,\VC\left(x^{2}\right)\right)}{\tau}+\sum_{k\geq2}\frac{Z\left(\VDc;\VC\left(x^{k}\right),\VC\left(x^{2k}\right)\right)}{k\,\VC\left(x^{k}\right)}
\right|_{x=\rho}\right).
\end{align*}
From $0<\tau<\delta$ and the fact that the derivative in $\frac{\partial}{\partial x}H(\rho,\tau)$
is a derivative of a formal power series with positive coefficients evaluated at $\rho>0$,  it follows that
both derivatives are strictly positive and hence do not vanish.

Finally, the aperiodicity of $\VC(x)$ is easily seen from the fact that $\vc_1\neq 0$ and $\vc_2\neq 0$.
\end{proof}

\begin{thm}[asymptotic number of connected outerplanar graphs]\label{thm:Casympt}
The function $\C\left(x\right)$
has a singular expansion of the form\begin{equation}
\C\left(x\right)=\C\left(\rho\right)+\sum_{k\geq2}\C_{k}X^{k},\qquad X:=\sqrt{1-\frac{x}{\rho}},\label{eq:Casympt}\end{equation}
 with constants $\C_{k}$, $k\geq2$, which can be computed from the
constants $\VC_{k}$, and with $\rho$ as in Lemma \ref{lem:VCsingularity}.
Hence,\[
\con_{n}\sim\frac{3\C_{3}}{4\sqrt{\pi}}n^{-5/2}\rho^{-n}.\]

\end{thm}
\begin{proof}
Recall Formula (\ref{eq:Cogf}) for the ordinary generating function
for connected outerplanar graphs,\[
\C\left(x\right)=\VC\left(x\right)+Z\left(\Dc;\VC\left(x\right)\right)-Z\left(\VDc;\VC\left(x\right)\right).\]
Since $\tau<\delta$, it is clear that the dominant singularity of
$\C\left(x\right)$ is the same as $\VC\left(x\right)$ \cite[Cha. VI.6]{FlajoletSedgewick}. The singular expansion of $\C\left(x\right)$ around $\rho$
can then be obtained by injecting the singular expansion of $\VC\left(x\right)$ into Formula (\ref{eq:Cogf}):
\begin{equation}
\begin{array}{rcl}
\C\left(x\right) & = & {\displaystyle \sum_{k\geq0}\VC_{k}X^{k}+Z\left(\Dc;\sum_{k\geq0}\VC_{k}X^{k},\VC\left(\rho^{2}\left(1-X^{2}\right)^{2}\right),\VC\left(\rho^{3}\left(1-X^{2}\right)^{3}\right),\ldots\right)}\\
&  & {\displaystyle -Z\left(\VDc;\sum_{k\geq0}\VC_{k}X^{k},\VC\left(\rho^{2}\left(1-X^{2}\right)^{2}\right),\VC\left(\rho^{3}\left(1-X^{2}\right)^{3}\right),\ldots\right)}.\end{array}\label{eq:Csubst}\end{equation}
Developing in terms of $X$ (around $X=0$) gives a singular
expansion $\C\left(x\right)=\sum_{k\geq0}\C_{k}X^{k}.$

It remains to check that $\C_{1}=0$ and $\C_3\neq 0$. From (\ref{eq:Csubst}) it
is clear that\[
\C_{1}=\VC_{1}+\VC_{1}\left.\frac{\partial}{\partial s_{1}}Z\left(\Dc\right)\right|_{\left(s_{1},s_{2}\right)=(\tau,\VC\left(\rho^{2}\right))}-\VC_{1}\left.\frac{\partial}{\partial s_{1}}Z\left(\VDc\right)\right|_{\left(s_{1},s_{2}\right)=(\tau,\VC\left(\rho^{2}\right))}.\]
From (\ref{eq:rooting}) we know $s_{1}\frac{\partial}{\partial s_{1}}Z\left(\Dc\right)=Z\left(\VDc\right),$
so that\[
\C_{1}=\VC_{1}\left.\left(1+\frac{Z\left(\VDc\right)}{s_{1}}-\frac{\partial}{\partial s_{1}}Z\left(\VDc\right)\right)\right|_{\left(s_{1},s_{2}\right)=(\tau,\VC\left(\rho^{2}\right))}.\]
On the other hand, Equation~(\ref{eq:HforVC}) implies that \begin{equation}
\frac{\partial}{\partial y}H\left(x,y\right)=\left(H\left(x,y\right)+y\right)\left(\frac{1}{y}\,\frac{\partial}{\partial s_{1}}Z\left(\VDc;y,\VC\left(x^{2}\right)\right)-\frac{1}{y^{2}}\, Z\left(\VDc;y,\VC\left(x^{2}\right)\right)\right)-1.\label{eq:Hdy}\end{equation}
By Equation~(\ref{eq:posimplicitfunc}) and Lemma \ref{lem:VCsingularity}, $0=\frac{\partial}{\partial y}H\left(\rho,\tau\right)=\left.\frac{\partial}{\partial s_{1}}Z\left(\VDc\right)-\frac{1}{s_{1}}Z\left(\VDc\right)-1\right|_{\left(s_{1},s_{2}\right)=(\tau,\VC\left(\rho^{2}\right))}$.
Thus, $\C_{1}=0$.
Assume $\C_3=0$. Then the expansion (\ref{eq:Casympt}) yields $\con_n \sim O(n^{-k/2-1}) \rho^{-n}$ for some odd number $k\geq 5$. This contradicts $n \con_n\geq\vc_n \sim - \VC_1/\left(2 \sqrt{\pi}\right) n^{-3/2}\rho^{-n}$
(by Lemma \ref{lem:VCsingularity}).
\end{proof}

\begin{thm}
[asymptotic number of outerplanar graphs]\label{thm:Out_asympt}The function $\Out\left(x\right)$
has a singular expansion of the form\[
\Out\left(x\right)=\Out\left(\rho\right)+\sum_{k\geq2}\Out_{k}X^{k},\qquad X:=\sqrt{1-\frac{x}{\rho}},\]
 where $\rho$ is as in Lemma \ref{lem:VCsingularity}, and where the constants 
$\Out_{k}$, $k\geq2$, can be computed from the constants $\C_{k}$, in particular $\Out_{3}=\Out\left(\rho\right)\C_{3}$. Furthermore, $\out_{n}$
has the asymptotic estimate\[
\out_{n}=\sum_{k\geq1}\binom{n+k-\frac{1}{2}}{n}\Out_{2k+1}\rho^{-n},\]
and in particular\[
\out_{n}\sim\frac{3\Out_{3}}{4\sqrt{\pi}}n^{-5/2}\rho^{-n}.\]

\end{thm}
\begin{proof}
Recall Formula (\ref{eq:Outogf}) for the ordinary generating function
for outerplanar graphs,\[
\Out\left(x\right)=\exp\left(\sum_{k\geq1}\frac{1}{k}\C\left(x^{k}\right)\right).\]
As the exponential function  $\exp( \cdot)$ is regular, the dominant singularity of $\Out\left(x\right)$ is the same as
$\C\left(x\right)$. Replacing $\C\left(x\right)$ by its singular
expansion (\ref{eq:Casympt}) and $x^{k}$ by $\rho^{k}\left(1-X^{2}\right)^{k}$
for $k\geq2$, we get\[
\Out\left(x\right)=\exp\left(\C\left(\rho\right)+\sum_{k\geq2}\C_{k}X^{k}+\sum_{k\geq2}\frac{1}{k}\C\left(\rho^{k}\left(1-X^{2}\right)^{k}\right)\right),\]
from which the singular expansion
of $\Out\left(x\right)$ can be computed. Then, by Lemma \ref{lem:basicscale} we
derive the asymptotic estimate of $\out_{n}$.
\end{proof}

Finally, using the same techniques as for the general case, we can
compute the asymptotic estimate of bipartite outerplanar graphs, given
in Theorem \ref{thm:bdisgrowthrate}.

\subsection{\label{sec:approxgrowth}Numerical approximation of the growth constants}

As far as we know, the computation of analytic expressions for growth constants
has not been possible for some classes of unlabeled structures that are even simpler than
outerplanar graphs, for example, for unembedded trees, see
\cite[Sec.~VII.2.3]{FlajoletSedgewick} and~\cite{Otter}.
Nevertheless, we can simplify the problem by reducing it to one variable, and
 provide numerical estimates of the growth constants.
With Formula (\ref{eq:Hdy}) for $\frac{\partial}{\partial y}H\left(\rho,\tau\right)$ and the explicit formula for $Z\left(\VDc\right)$ from Corollary \ref{cor:VD} the equation $\frac{\partial}{\partial y}H\left(\rho,\tau\right)=0$
becomes \begin{equation}
\tau\left(1+\VC\left(\rho^{2}\right)\left(\VC\left(\rho^{2}\right)-3\right)-\frac{\VC\left(\rho^{2}\right)^{2}\left(\tau-3\right)}{\sqrt{\tau^{2}-6\tau+1}}-\sqrt{\VC\left(\rho^{2}\right)^{2}-6\VC\left(\rho^{2}\right)+1}\right)=8\VC\left(\rho^{2}\right)^{2}.\label{eq:Hdysimple}\end{equation}
With algebraic elimination~\cite[App. B.1]{FlajoletSedgewick},
Equation~(\ref{eq:Hdysimple}) can be reformulated as a system
of polynomial equations, regarding $\VC\left(\rho^{2}\right)$ as
a fixed value. We obtain a polynomial equation of degree 8 in $\tau$
with coefficients $p_{i}\left(\rho\right)$
(depending on $\VC\left(\rho^{2}\right)$), $i=0,\ldots,8$,
\begin{equation}
p_{0}\left(\rho\right)+p_{1}\left(\rho\right)\tau+p_{2}\left(\rho\right)\tau^{2}+p_{3}\left(\rho\right)\tau^{3}+p_{4}\left(\rho\right)\tau^{4}+p_{5}\left(\rho\right)\tau^{5}+p_{6}\left(\rho\right)\tau^{6}+p_{7}\left(\rho\right)\tau^{7}+p_{8}\left(\rho\right)\tau^{8}=0.\label{eq:taupolynom}\end{equation}
The solutions of (\ref{eq:taupolynom}) do not need to satisfy
Equation~(\ref{eq:Hdysimple}), but every $\tau$ that is a solution
of (\ref{eq:Hdysimple}) is also a solution of (\ref{eq:taupolynom}) (see \cite{Vigerske} for the details). We denote the solutions of (\ref{eq:taupolynom})
by $\tau_{1}\left(\rho\right),\ldots,\tau_{8}\left(\rho\right)$.
It remains to solve the equations\[
H\left(\rho,\tau_{i}\left(\rho\right)\right)=0,\qquad i=1,\ldots,8,\]
 and to pick the correct solution $\rho$. Since $H\left(x,y\right)$
depends on $\VC\left(x\right)$, which we do not know explicitly,
and since it includes also an infinite sum that we were not able to
simplify, we can only approximate the solutions of $H\left(\rho,\tau_{i}\left(\rho\right)\right)=0$
by truncating the infinite sum in $H\left(x,y\right)$ at some
index $m$ and replacing $\VC\left(x\right)$ with $\VC^{\left[m\right]}\left(x\right):=\sum_{n=1}^{m}\vc_{n}x^{n}$
for known coefficients $\vc_{1}$, $\ldots$, $\vc_{m}$. That is,
we search for roots of the functions\[
\tilde{H}_{i}^{\left[m\right]}\left(\rho\right):=\rho\exp\left(\frac{Z\left(\VDc;\tau_{i}\left(\rho\right),\VC^{\left[m\right]}\left(\rho^{2}\right)\right)}{\tau_{i}\left(\rho\right)}+\sum_{k=2}^{m}\frac{Z\left(\VDc;\VC^{\left[m\right]}\left(\rho^{k}\right),\VC^{\left[m\right]}\left(\rho^{2k}\right)\right)}{k\:\VC^{\left[m\right]}\left(\rho^{k}\right)}\right)-\tau_{i}\left(\rho\right),\]
 $i=1,\ldots,8$, in the interval $\left(0,1\right)$.
We solve the equation $\tilde{H}_{i}^{\left[m\right]}\left(\rho\right)=0$
for $m=25$ numerically, select the correct root, and obtain the estimates
\[
\rho\approx0.1332694\quad\textrm{and}\quad\tau\approx0.1707560.\]
The residuals in the equations $\tilde{H}_{i}^{\left[m\right]}\left(\rho,\tau\right)=0$
and $\frac{\partial}{\partial y}\tilde{H}_{i}^{\left[m\right]}\left(\rho,\tau\right)=0$
have an order of $10^{-58}$.
Table \ref{cap:H3graph} shows approximations of $\rho$ for several values of $m$.

\begin{table}[!ht]
\begin{center}\begin{tabular}{|c|l|}
\hline
$m$&
\multicolumn{1}{c|}{approximation of $\rho$}\tabularnewline
\hline
\hline
$1$&
\textbf{0.13}461876886110181369...\tabularnewline
\hline
$4$&
\textbf{0.1332}7064317786556821...\tabularnewline
\hline
$8$&
\textbf{0.133269432}88029243729...\tabularnewline
\hline
$16$&
\textbf{0.1332694326674468}2071...\tabularnewline
\hline
$25$&
\textbf{0.13326943266744680944}...\tabularnewline
\hline
\end{tabular}\end{center}
\caption{\label{cap:H3graph}The accuracy is improved by increasing the order of truncation.}
\end{table}

We can now estimate the coefficients in the singular expansions
of $\VC\left(x\right)$, $\C\left(x\right)$ and $\Out\left(x\right)$.
In particular $\VC_1\approx -0.0255905$, $C_3\approx 0.0179720$ and $G_3\approx 0.0215044$.

The growth constant for bipartite outerplanar graphs can also be
estimated in the same way as $\rho$, and we get 
$\rho_{\mathrm{\it{b}}}\approx0.218475$ (see \cite{Vigerske} for details).

\section{\label{sec:randomgraphs}Random unlabeled outerplanar graphs}

This section investigates typical properties of a random unlabeled outerplanar
graph with $n$ vertices. We first discuss the probability of being connected, and the number and type of components, 
and then proceed with the distribution
of the number of edges.

\subsection{Connectedness, components, and isolated vertices}
We start with the proof of Theorem {\ref{thm:typprop}} (1) on the probability that a random outerplanar graph is connected.

\begin{proof}
[Proof of Theorem {\ref{thm:typprop}} (1)]The probability that
a random outerplanar graph on $n$ vertices is connected is exactly $\con_{n}/\out_{n}$.
The asymptotic estimates for $\con_{n}$ and $\out_{n}$ from Theorem~\ref{thm:Casympt} and Theorem~\ref{thm:Out_asympt} yield $\con_{n}/\out_{n}\sim\C_{3}/\Out_{3}\approx0.845721$.
\end{proof}
The number of components can be studied by augmenting the generating
function for outerplanar graphs with a variable that counts the number of components.

\begin{proof}
[Proof of Theorem {\ref{thm:typprop}} (2)]Let $\kappa_{n}$
denote the number of components in a random outerplanar graph on $n$
vertices and let $\Out\left(x,u\right):=\exp\left(\sum_{k\geq1}\frac{1}{k}u^{k}\C\left(x^{k}\right)\right)$ be
the generating function for outerplanar graphs, where  the additional variable $u$ marks
the number of components. Thus,
the pro\-ba\-bi\-li\-ty that an outerplanar graph has $k$ components
is $\mathbb{P}\left[\kappa_{n}=k\right]=\left[x^{n}u^{k}\right]\Out\left(x,u\right)/\out_{n},$
and the expected number of components is \[
\mathbb{E}\left[\kappa_{n}\right]=\frac{1}{\out_{n}}\sum_{k\geq1}k\left[x^{n}u^{k}\right]\Out\left(x,u\right)=\frac{1}{\out_{n}}\left[x^{n}\right]\frac{\partial}{\partial u}\Out\left(x,1\right)=\frac{1}{\out_{n}}\left[x^{n}\right]\Out\left(x\right)\sum_{k\geq1}\C\left(x^{k}\right).\]
By asymptotic expansion around $x=\rho$, we obtain \[
\left[x^{n}\right]\Out\left(x\right)\sum_{k\geq1}\C\left(x^{k}\right)\sim\Out\left(\rho\right)\C_{3}\left(1+\sum_{r\geq1}\C\left(\rho^{r}\right)\right)\frac{1}{\Gamma\left(-3/2\right)}n^{-5/2}\rho^{-n},\]
 which together with Theorem \ref{thm:Out_asympt}, more precisely $\out_{n}\sim\Out\left(\rho\right)\C_{3}\frac{1}{\Gamma\left(-3/2\right)}n^{-5/2}\rho^{-n},$
yields \[
\mathbb{E}\left[\kappa_{n}\right]\sim1+\sum_{r\geq1}\C\left(\rho^{r}\right)\approx1.17847.\qedhere\]
\end{proof}
Given a family $\mathcal{A}$ of connected outerplanar graphs, we
can make the following statements about the probability that a random
outerplanar graph has exactly $k$ components in $\mathcal{A}$. Denote
the number of graphs in $\mathcal{A}$ that have exactly $n$ vertices
by $a_{n}$, and let $A\left(x\right):=\sum_{n}a_{n}x^{n}$.

\begin{thm}
\label{thm:nrcomponentsspecifictype}Given an outerplanar graph $G$ with $n$ vertices, let $\kappa_{n}^{\mathcal{A}}$
be the number of connected components of $G$
belonging to $\mathcal{A}$. If the radius of convergence $\alpha$
of $A\left(x\right)$ is strictly larger than $\rho$,
that is, $a_n$ is exponentially smaller than $c_n$, then
the probability that a random outerplanar graph with $n$ vertices has exactly $k\geq0$
components belonging to $\mathcal{A}$ converges to a discrete law:\[
\mathbb{P}\left[\kappa_{n}^{\mathcal{A}}=k\right]\sim Z\left(S_{k};A\left(\rho\right)\right)\exp\left(-\sum_{r\geq1}\frac{1}{r}A\left(\rho^{r}\right)\right)\; ,\]
and the expected number of components belonging to $\mathcal{A}$
in a random outerplanar graph with $n$ vertices is \[
\mathcal{\mathbb{E}}\left[\kappa_{n}^{\mathcal{A}}\right]\sim\sum_{r\geq1}A\left(\rho^{r}\right).\]

\end{thm}
\begin{proof}
Let $G^\mathcal{A}(x,u)$ be the generating function for outerplanar graphs, where the additional variable $u$ marks the number of components belonging to $\mathcal{A}$,
\begin{eqnarray*}
\Out^{\mathcal{A}}\left(x,u\right) & := & \exp\left(\sum_{k\geq1}\frac{1}{k}\left(u^{k}A\left(x^{k}\right)+\left(\C\left(x^{k}\right)-A\left(x^{k}\right)\right)\right)\right)\\
& = & \Out\left(x\right)\exp\left(\sum_{k\geq1}\frac{u^{k}-1}{k}A\left(x^{k}\right)\right).
\end{eqnarray*}
Hence, $\mathbb{P}\left[\kappa_{n}^{\mathcal{A}}=k\right]=\left[x^{n}u^{k}\right]\Out^{\mathcal{A}}\left(x,u\right)/\out_{n}$.
Since $A\left(x\right)$ is analytic at $\rho$, the dominant singularity
of $\Out^{\mathcal{A}}\left(x,u\right)$ for fixed $u$ is determined
by $\Out\left(x\right)$. Thus,
\[
[x^nu^k]\Out^{\mathcal{A}}\left(x,u\right)\mathop{\sim}_{n\to \infty}
[u^k]\exp\left(\sum_{k\geq 1}\frac{u^k-1}{k}A(\rho^k)\right)[x^n]G(x),
\] i.e., \[
\mathbb{P}\left[\kappa_{n}^{\mathcal{A}}=k\right]\mathop{\sim}_{n\to\infty}\left[u^{k}\right]\exp\left(\sum_{k\geq1}\frac{u^{k}-1}{k}A\left(\rho^{k}\right)\right)=Z\left(S_{k};A\left(\rho\right)\right)\exp\left(-\sum_{k\geq1}\frac{1}{k}A\left(\rho^{k}\right)\right).\]

For the expectation of $\kappa_{n}^{\mathcal{A}}$ we again use \[
\mathbb{E}\left[\kappa_{n}^{\mathcal{A}}\right]=\frac{1}{\out_{n}}\left[x^{n}\right]\frac{\partial}{\partial u}\Out^{\mathcal{A}}\left(x,1\right)=\frac{1}{\out_{n}}\left[x^{n}\right]\Out\left(x\right)\sum_{k\geq1}A\left(x^{k}\right).\]
The statement follows from the analyticity of $A\left(x\right)$ at $\rho$ and Theorem~\ref{thm:Out_asympt}.
\end{proof}
The asymptotic distribution of the number of isolated vertices in
a random outerplanar graph can now be easily computed, as stated in Theorem \ref{thm:typprop} (3).

\begin{proof}[Proof of Theorem {\ref{thm:typprop}} (3)]
Let $\mathcal{A}$
be the family consisting of the graph that is a single vertex, i.e., $A\left(x\right)=x$.
By Theorem \ref{thm:nrcomponentsspecifictype}, $\mathbb{P}\left[\kappa_{n}^{\mathcal{A}}=k\right]\sim\rho^{k}/\left(1-\rho\right)$,
since $Z\left(S_{k};A\left(\rho\right)\right)=\rho^{k}$ and $\sum_{r}\frac{1}{r}A\left(\rho^{r}\right)=\log\left(1-\rho\right)$.
Hence, the distribution of the number of isolated vertices
$\kappa_{n}^{\mathcal{A}}$ is asymptotically a geometric law with parameter $\rho$.
\end{proof}
Other consequences of Theorem \ref{thm:nrcomponentsspecifictype}
concern the number of two-connected components and the number of bipartite
components in a random outerplanar graph.

\begin{cor}[two-connected components]
In a random outerplanar graph, the expected number of connected components that are two-connected is asymptotically $\sum_{k\geq1}\D\left(\rho^{k}\right)\approx0.175054.$
\end{cor}
\begin{proof}
Let $\mathcal{A}:=\Dc$ be the family of dissections, $A\left(x\right)=\D\left(x\right)$.
The radius of convergence of $\D\left(x\right)$ is $\delta>\rho$ (Lemma \ref{lem:VCsingularity}).
Hence, by Theorem \ref{thm:nrcomponentsspecifictype}, $\mathbb{E}\left[\kappa_{n}^{\Dc}\right]=\sum_{k\geq1}\D\left(\rho^{k}\right).$
\end{proof}
\begin{cor}
[number of bipartite components]In a random outerplanar graph, the expected number of connected components that are bipartite
 is asymptotically
$\sum_{k\geq1}\bC\left(\rho^{k}\right)\approx0.175427$, where $\bC\left(x\right)$ is the generating function for bipartite connected outerplanar graphs.
\end{cor}
\begin{proof}
We apply Theorem \ref{thm:nrcomponentsspecifictype} with $\mathcal{A}=\bCc$.
\end{proof}

\subsection{Number of edges}
In this section, we analyze the distribution of the number of edges in a random
outerplanar graph. 
To do this, we add a variable $y$ whose power (in the cycle index sums
and generating functions) indicates the number of edges.
For a graph $G$ on $n$ vertices and $m$ edges, and with the automorphism
group $\Gamma\left(G\right)$ (acting on the vertices), we define\[
Z\left(G;s_{1},s_{2},\ldots;y\right):=Z\left(\Gamma\left(G\right);s_{1},s_{2},\ldots;y\right):=y^{m}\frac{1}{\left|\Gamma\left(G\right)\right|}\sum_{\alpha\in\Gamma\left(G\right)}\prod_{k=1}^{n}s_{k}^{j_{k}\left(\alpha\right)}.\]
Taking the number of edges into account in the calculations of
Section~\ref{sec:cycleindexsums}, the cycle index sums for all encountered
families of outerplanar graphs can be derived with the additional
variable $y$ (see \cite{Vigerske} for more details).

{\allowdisplaybreaks\begin{align*}
Z\left(\OEDc\right) & =\frac{s_{1}}{2\left(1+y\right)}\left(s_{1}y+1-\sqrt{\left(s_{1}y-1\right)^{2}-4s_{1}y^{2}}\right),\\
Z^{+}\left(\REDc\right) & =s_{1}\frac{\left(s_{1}y-1\right)\left(s_{1}^{2}y^{2}-2s_{1}y^{2}-1\right)-\left(1+s_{1}y\right)\sqrt{\left(s_{1}^{2}y^{2}-1\right)^{2}-4s_{1}^{2}y^{4}}}{2\left(s_{1}^{2}y^{3}+s_{1}^{2}y^{2}+y-1\right)},\\
Z^{-}\left(\REDc\right) & =\frac{s_{2}\left(2y^{2}+s_{2}y^{3}-y\right)-s_{1}\left(s_{2}y^{2}+2s_{2}y^{3}-1\right)-\left(s_{1}+s_{2}y\right)\sqrt{\left(s_{2}y^{2}-1\right)^{2}-4s_{2}y^{4}}}{2\left(y+s_{2}y^{2}+s_{2}y^{3}-1\right)},\\
Z\left(\REDc\right) & =\frac{1}{2}\left(Z^{+}\left(\REDc\right)+Z^{-}\left(\REDc\right)\right),\\
Z\left(\IDc\right) & =\frac{y}{4}\left(\left(\frac{Z\left(\OEDc\right)}{s_{1}y}-s_{1}\right)^{2}+\left(\frac{Z\left(\OEDc;s_{2};y^{2}\right)}{s_{2}^{2}y^{2}}-1\right)\left(s_{1}^{2}+s_{2}\right)+\left(\frac{Z^{-}\left(\REDc\right)}{s_{2}y}-1\right)^{2}s_{2}\right),\\
Z\left(\SDc\right) & =\frac{Z\left(\OEDc;s_{1}^{2};y^{2}\right)}{2s_{1}^{2}y}+\frac{s_{1}^{2}+s_{2}}{4s_{2}^{2}y}Z\left(\OEDc;s_{2};y^{2}\right)-\frac{Z^{+}\left(\REDc;s_{1}^{2};y^{2}\right)}{4s_{1}^{2}y}+\frac{Z^{-}\left(\REDc;s_{1}^{2},s_{2}^{2};y^{2}\right)}{4s_{2}y}-\frac{s_{1}^{2}+s_{2}}{2}y,\\
Z\left(\OFDc\right) & =-\left\{ \sum_{d}\frac{\varphi\left(d\right)}{d}\log\left(1-\frac{Z\left(\OEDc;s_{d};y^{d}\right)}{s_{d}}\right)\right\} -\frac{Z\left(\OEDc\right)}{s_{1}}-\frac{1}{2}\left(\left(\frac{Z\left(\OEDc\right)}{s_{1}}\right)^{2}+\frac{Z\left(\OEDc;s_{2};y^{2}\right)}{s_{2}}\right),\\
Z\left(\FDc\right) & =\frac{1}{2}Z\left(\OFDc\right)+\frac{1}{2s_{2}}\left(s_{1}Z^{-}\left(\REDc\right)+\frac{s_{1}^{2}}{2s_{2}}Z\left(\OEDc;s_{2};y^{2}\right)+\frac{1}{2}Z^{-}\left(\REDc\right)^{2}\right)\:\frac{Z\left(\OEDc;s_{2};y^{2}\right)}{s_{2}-Z\left(\OEDc;s_{2};y^{2}\right)},\\
Z\left(\SFDc\right) & =\frac{1}{s_{1}^{2}y}Z\left(\OEDc;s_{1}^{2};y^{2}\right)-\frac{1}{2s_{1}^{2}y}Z^{+}\left(\REDc;s_{1}^{2},s_{2}^{2};y^{2}\right)+\frac{1}{2s_{2}y}Z^{-}\left(\REDc;s_{1}^{2},s_{2}^{2};y^{2}\right)-\frac{y}{2}\left(s_{1}^{2}+s_{2}\right),\\
Z\left(\Dc\right) & =\frac{1}{2}\left(s_{1}^{2}+s_{2}\right)y+Z\left(\FDc\right)-Z\left(\SFDc\right)-Z\left(\IDc\right)+2Z\left(\SDc\right),\\
Z\left(\VDc\right) & =s_{1}\frac{\partial}{\partial s_{1}}Z\left(\Dc;s_{1},s_{2},\ldots;y\right),\\
Z(\VCc) & =s_{1}\exp\left(\sum_{k\geq1}\frac{1}{k}\frac{Z\left(\VDc;Z\left(\VCc;s_{k},s_{2k};y^{k}\right),Z\left(\VCc;s_{2k},s_{4k};y^{2k}\right),\ldots;y^{k}\right)}{Z\left(\VCc;s_{k},s_{2k};y^{k}\right)}\right),\\
Z\left(\Cc\right) & =Z(\VCc)+Z\left(\Dc;Z(\VCc)\right)-Z\left(\VDc;Z(\VCc)\right),\\
Z\left(\Outc\right) & =\exp\left(\sum_{k\geq1}\frac{Z\left(\Cc;s_{k},s_{2k},\ldots;y^{k}\right)}{k}\right).\end{align*}
}%
Similarly as in Section~\ref{sec:cycleindexsums}, the coefficients counting outerplanar graphs with respect to the number of vertices and the number of edges can be extracted in polynomial time from the expressions of the cycle index sums, see~\cite{Vigerske} for a table.

With the help of Theorem \ref{thm:centrallaw}, we
can prove Theorem \ref{thm:edges}
giving the limit distributions of the number of edges in a random
dissection and in a random outerplanar graph, respectively.

\begin{proof}[Proof of Theorem {\ref{thm:edges}}]
We start with the limit distribution of the number of edges in a two-connected outerplanar graph.
The generating function for
oriented outer-edge rooted dissections that additionally counts edges
is \[
\OED\left(x,y\right)=\frac{x}{2\left(y+1\right)}\left(xy+1-\sqrt{\left(xy-1\right)^{2}-4xy^{2}}\right).\]
The singularities of $\OED\left(x,y\right)$ are determined by the
equation $\left(xy-1\right)^{2}-4xy^{2}=0$. Hence, for $y$ close
to 1, the dominant singularity of $x\mapsto \OED\left(x,y\right)$ is at
$\delta\left(y\right)=2+1/y-2\sqrt{1+1/y}.$
With the
same arguments as before, $\delta\left(y\right)$ is also the dominant
singularity of the generating functions for vertex rooted and unrooted
dissections. Furthermore, $\delta'\left(1\right)=-1+\sqrt{2}/2$ and
$-\frac{\delta''\left(1\right)}{\delta\left(1\right)}-\frac{\delta'\left(1\right)}{\delta\left(1\right)}+\left(\frac{\delta'\left(1\right)}{\delta\left(1\right)}\right)^{2}=\sqrt{2}/8\neq0,$
so that the variance condition (in Theorem \ref{thm:centrallaw}) on $\delta\left(y\right)$ is satisfied.
Hence, Theorem \ref{thm:centrallaw} yields the statement for dissections.

We now determine the distribution of the number of edges in a rooted
connected outerplanar graph. The generating function $\VC\left(x,y\right)$
is implicitly defined by\[
\VC\left(x,y\right)=x\exp\left(\sum_{k\geq1}\frac{Z\left(\VDc;\VC\left(x^{k},y^{k}\right);y^{k}\right)}{k\,\VC\left(x^{k},y^{k}\right)}\right).\]
In order to apply the singular implicit functions theorem~\ref{thm:positiveimplicitfunctions} for the function $x\mapsto\VC\left(x,y\right)$ with a fixed $y$ close to 1, we define 
\[
H\left(x,y,z\right):=x\exp\left(\frac{Z\left(\VDc;z,\VC\left(x^{2},y^{2}\right);y\right)}{z}+\sum_{k\geq2}\frac{Z\left(\VDc;\VC\left(x^{k},y^{k}\right),\VC\left(x^{2k},y^{2k}\right);y^{k}\right)}{k\:\VC\left(x^{k},y^{k}\right)}\right)-z\]
 and search for a solution $\left(x,z\right)=\left(\rho\left(y\right),\tau\left(y\right)\right)$
of the system\begin{equation}
H\left(x,y,z\right)=0,\qquad\frac{\partial}{\partial z}H\left(x,y,z\right)=0,\label{eq:Hsystemedges}\end{equation}
such that $\left(\rho\left(y\right),\tau\left(y\right)\right)$ is in the analyticity domain of $(x,z)\mapsto H(x,y,z)$.

For $y=1$, the solution is at $x=\rho$, $z=\tau$ by Lemma \ref{lem:VCsingularity}. Then the classical implicit functions theorem, applied to the system~(\ref{eq:Hsystemedges}), ensures that the solution $(\rho,1,\tau)$ can be extended into solutions $(\rho(y),y,\tau(y))$ for $y$ close to 1, where the functions $\rho(y)$ and $\tau(y)$ are analytic in a neighbourhood of 1.
To apply the classical implicit function theorem  on system (\ref{eq:Hsystemedges}), it remains to check that the determinant of the Jacobian
of system (\ref{eq:Hsystemedges}), with respect to $x$ and $z$,\[
\left(\begin{array}{cc}
{\displaystyle \frac{\partial}{\partial x}H\left(x,y,z\right)} & {\displaystyle \frac{\partial}{\partial z}H\left(x,y,z\right)}\medskip\\
{\displaystyle \frac{\partial}{\partial x}\frac{\partial}{\partial z}H\left(x,y,z\right)} & {\displaystyle \frac{\partial}{\partial z}\frac{\partial}{\partial z}H\left(x,y,z\right)}\end{array}\right),\]
 does not vanish at $\left(x,y,z\right)=\left(\rho\left(1\right),1,\tau\left(1\right)\right)$.
This is clear, since from Lemma \ref{lem:VCsingularity} we have
$\frac{\partial}{\partial z}H\left(\rho\left(1\right),1,\tau\left(1\right)\right)=0$,
$\frac{\partial}{\partial x}H\left(\rho\left(1\right),1,\tau\left(1\right)\right)\neq0$,
and $\frac{\partial^{2}}{\partial z^{2}}H\left(\rho\left(1\right),1,\tau\left(1\right)\right)\neq0$.
Hence, there exist analytic functions $\rho\left(y\right)$ and $\tau\left(y\right)$
such that\begin{equation}
H\left(\rho\left(y\right),y,\tau\left(y\right)\right)=0,\qquad\frac{\partial}{\partial z}H\left(\rho\left(y\right),y,\tau\left(y\right)\right)=0,\label{eq:Hsystemedgessol}\end{equation}
$\frac{\partial^{2}}{\partial z^{2}}H\left(\rho\left(y\right),y,\tau\left(y\right)\right)\neq0$,
and $\frac{\partial}{\partial x}H\left(\rho\left(y\right),y,\tau\left(y\right)\right)\neq0$
for $y$ close to one. In addition, these solutions are in the analyticity domain of $(x,z)\mapsto H(x,y,z)$ for $y$ close to 1, by analyticity of $(x,y,z)\mapsto H(x,y,z)$ at $(\rho,1,\tau)$.
Next, the singular implicit functions theorem~\ref{thm:positiveimplicitfunctions}
yields a singular expansion $\VC\left(x,y\right)=\sum_{k\geq0}\VC_{k}\left(y\right)(\sqrt{1-x/\rho\left(y\right)})^{k}$
with coefficients $\VC_{k}\left(y\right)$ analytic at $y=1$ and verifying 
$\VC_{1}\left(y\right)\neq0$ for $y$ close to 1.

To find $\rho'\left(1\right)$ and $\rho''\left(1\right)$ we compute the first and second derivatives of the equations in (\ref{eq:Hsystemedgessol})
with respect to $y$, and express $\rho'\left(y\right)$ and $\rho''\left(y\right)$ in terms of $\rho\left(y\right)$, $\tau\left(y\right)$, and the partial
derivatives of $H\left(x,y,z\right)$ at $\left(x,z\right)=\left(\rho\left(y\right),\tau\left(y\right)\right)$.
Using the approximated values from Section \ref{sec:approxgrowth} we obtain $\rho'\left(1\right)\approx-0.206426$,
$\rho''\left(1\right)\approx0.495849$, and $-\frac{\rho''\left(1\right)}{\rho\left(1\right)}-\frac{\rho'\left(1\right)}{\rho\left(1\right)}+(\frac{\rho'\left(1\right)}{\rho\left(1\right)})^{2}\approx0.227504\neq0$.
Theorem \ref{thm:centrallaw} implies that the distribution of the number of edges
in a random rooted connected outerplanar graph with $n$ vertices 
asymptotically follows a Gaussian law with mean $\mu n$ and variance $\sigma^2n$, where $\mu=-\frac{\rho'\left(1\right)}{\rho\left(1\right)}\approx1.54894$
and $\sigma^2\approx0.227504.$ The same holds for unrooted connected outerplanar graphs
and for outerplanar graphs, since their generating
functions have the same dominant singularity.
\end{proof}

\section{Conclusion}
A summary of the estimated growth constants and other parameters for
unlabeled outerplanar graphs is presented in Table \ref{cap:summary}.
For comparison we also include the corresponding labeled quantities derived in \cite{BodirskyGimenezKangNoy}. 
Observe that in the two-connected case the estimated
quantities for the unlabeled and labeled structures do not differ,
since their dominant singularity is determined by the same equation
(compare Theorem \ref{thm:twoconn} and the formula for $B\left(x\right)$
in \cite{BodirskyGimenezKangNoy}). 
\begin{table}[!ht]
\begin{center}\begin{tabular}{|c||c|c|c|}
\hline 
&
dissections&
\multicolumn{2}{c|}{outerplanar graphs}\tabularnewline
&
unlabeled or labeled&
unlabeled&
labeled \tabularnewline
\hline
\hline 
growth constant&
$\delta^{-1}\approx5.82843$&
$\rho^{-1}\approx7.50360$&
$\lambda^{-1}\approx7.32098$\tabularnewline
\hline 
$\mathbb{P}\left[\textrm{connectivity}\right]$&
1&
$0.845721$&
$0.861667$\tabularnewline
$\mathbb{E}\left[\textrm{nr. of components}\right]$&
1&
$1.17847$&
-\tabularnewline
\hline 
distr. of nr. of isolated vertices&
Dirac&
Geom$\left(\rho\right)$&
Poiss$\left(\lambda\right)$\tabularnewline
$\mathbb{E}\left[\textrm{nr. of isolated vertices}\right]$&
0&
$0.153761$&
$0.136593$\tabularnewline
\hline 
chromatic number&
3&
3&
3\tabularnewline
\hline
distr. of nr. of edges&
Gaussian&
Gaussian&
Gaussian\tabularnewline
$\mathbb{E}\left[\textrm{nr. of edges}\right]$&
$1.70711n$&
$1.54894n$&
$1.56251n$\tabularnewline
$\mathbb{V}\left[\textrm{nr. of edges}\right]$&
$0.176777n$&
$0.227504n$&
$0.223992n$\tabularnewline
\hline
\end{tabular}\end{center}

\caption{\label{cap:summary}Summary of growth constants, typical properties,
and limit laws for unlabeled and labeled dissections and outerplanar
graphs.}
\end{table}

\end{document}